\documentclass[11pt,a4paper]{article}
\ifdefined\pdfminorversion\pdfminorversion=4\fi
\ifdefined\pdfobjcompresslevel\pdfobjcompresslevel=0\fi
\usepackage[margin=1in]{geometry}
\usepackage[authoryear,round]{natbib}
\usepackage{amsmath,amssymb,mathtools,bm,mathrsfs,booktabs,array,enumitem,url,graphicx,amsthm}
\usepackage[hidelinks]{hyperref}
\setlist[itemize]{leftmargin=*,topsep=2pt,itemsep=1pt}
\newtheorem{theorem}{Theorem}
\newtheorem{lemma}[theorem]{Lemma}
\newtheorem{proposition}[theorem]{Proposition}
\newtheorem{corollary}[theorem]{Corollary}
\newtheorem{remark}[theorem]{Remark}
\newtheorem{definition}{Definition}
\newcounter{algorithm}
\newtheorem{localpolicy}[algorithm]{Algorithm}
\DeclareMathOperator{\Cov}{Cov}
\DeclareMathOperator{\diag}{diag}
\DeclareMathOperator{\vecop}{vec}
\DeclareMathOperator{\col}{col}
\newcommand{\E}{\mathbb E}
\newcommand{\Prob}{\mathbb P}
\newcommand{\K}{\mathcal K}
\newcommand{\A}{\mathcal A}
\newcommand{\F}{\mathcal F}
\newcommand{\Wien}{\mathcal W}
\newcommand{\id}{\mathrm I}
\newcommand{\Reg}{\mathcal R}
\newcommand{\norm}[1]{\left\lVert#1\right\rVert}
\newcommand{\tr}{\operatorname{tr}}

\let\standardappendix\appendix
\renewcommand{\appendix}{%
  \standardappendix
  \numberwithin{equation}{section}%
  \numberwithin{figure}{section}%
  \numberwithin{table}{section}%
}
\title{Expected Minimax Regret Rates for Adaptive LQG Control on Fixed Parameter Sets}
\author{Zhaobo Liu\thanks{Corresponding author. Email: \texttt{liuzhaobo@szu.edu.cn}.}\\[0.4em]
  {\small Institute for Advanced Study, Shenzhen University}\\
  {\small Shenzhen 518060, China}}
\date{1 November 2026}
\hypersetup{
  pdftitle={Expected Minimax Regret Rates for Adaptive LQG Control on Fixed Parameter Sets},
  pdfauthor={Zhaobo Liu}
}
\begin{document}
\setlength{\jot}{1pt}
\maketitle
\begin{abstract}
We study how a fixed model family determines expected minimax regret
in adaptive LQG control. The family is compact and affine, with partial
observation and unknown noise covariances. We quantify how information from residual covariances and input responses
can vanish relative to differences between optimal control laws. This relation yields
matching bounds that classify regret as bounded, logarithmic, sublinear
polynomial or linear. Polynomial exponents are rational.
For real algebraic model and cost data, a terminating procedure determines
the regret class and the exact exponent of polynomial growth before data
collection. A causal policy attains these rates using model estimation and bounded
probes, reselecting model sets from observations without a supplied
common stabilizer. For costs with fixed positive definite output and input
weights, polynomial exponents are at most one half.
\end{abstract}
{\small\noindent\textit{Keywords:} Adaptive LQG, minimax regret, partial observation, bounded probing, supervisory control\par}

\section{Introduction}\label{sec:intro}

Adaptive linear--quadratic--Gaussian (LQG) control seeks nearly optimal performance when the dynamics are unknown and only noisy outputs are observed \citep{campi1998,prandini2001}. Each input affects both the control cost and the information available for future decisions. Classical self-tuning regulators combine identification and feedback \citep{astrom1973}, with stochastic convergence and asymptotic performance guarantees established \citep{goodwin1981,kumar1983,chen1987,makila1984}. These asymptotic guarantees do not quantify the additional control cost incurred during learning over a given operating horizon. Finite-horizon regret measures the additional expected cost relative to an optimal controller that knows the model but has access to the same type of observations \citep{lale2020c,lale2020a}.

For fully observed linear--quadratic regulation (LQR), regret bounds have been established \citep{abbasi11,abeille17,dean18,cohen19}, with square-root regret attainable up to logarithmic factors for the general class \citep{simchowitz20}. Under partial observation, ExpCommit gives $\widetilde O(T^{2/3})$ regret and LQGOpt gives $\widetilde O(\sqrt T)$ regret under their respective assumptions \citep{lale2020c,lale2020a}. Certainty equivalence applies the optimal controller for an estimated model. Its LQG performance gap admits a local upper bound proportional to the squared model error \citep{mania19}, and additive exploration yields $\widetilde O(\sqrt T)$ regret \citep{athrey2024}. This raises the question: under what conditions can regret grow more slowly than the square root of the horizon?

Structural information can improve on these general rates. With a known stabilizing controller, logarithmic expected regret is attainable in LQR when the input matrix is known, or when the state matrix is known and the optimal gain has full row rank \citep[Theorems~2--3]{jedra2022}. Under an identifiability condition on structural side information, input perturbations yield polylogarithmic high-probability bounds for regret measured against an oracle trajectory \citep[Theorem~4]{shirani2020}. AdaptOn obtains polylogarithmic regret for strongly convex costs relative to the best controller in a specified persistently exciting class \citep{lale2020b}. Related work classifies regret on shrinking $T^{-1/4}$ neighborhoods in fully observed structured LQR and gives logarithmic guarantees on sufficiently small fixed neighborhoods under additional conditions \citep{liu2026}.

To assess how far structural information can reduce regret, regret lower bounds are also needed. For memoryless linear quadratic control, a singularity condition on Fisher information rules out logarithmic regret and yields a square-root lower bound \citep{ziemann2021a}. For LQR with a known state matrix and an affinely parameterized input matrix, an uninformative optimal policy yields a square-root local asymptotic minimax lower bound on shrinking parameter neighborhoods \citep[Theorem~5]{ziemann2021b}. Here, uninformativeness means that Fisher information under the optimal policy is singular along parameter directions in which the optimal controller changes to first order. Square-root lower bounds have also been established for a class of partially observed systems, with explicit dependence on system parameters \citep{ziemann2025}. Thus, structural information permits smaller regret in some settings, but square-root regret remains unavoidable in others even when such structure is used.

We study expected minimax regret on a given compact affine family of partially observed systems. We examine the information needed to distinguish models with different optimal control laws and the cost of obtaining that information. We relate differences between optimal laws to the performance loss caused by selecting the wrong law, and use this relation to characterize regret growth on the family.

Our main contributions are as follows.

1) We characterize expected minimax regret on the given model family, determining when it is bounded, logarithmic, sublinear polynomial, or linear, and identify the sharp rational exponent governing sublinear polynomial growth. For the logarithmic and sublinear polynomial regimes, we construct learning policies that attain these rates and establish matching lower bounds for arbitrary causal measurable randomized policies.

2) When the optimal law varies across the family, we establish regret guarantees over the entire family without prior knowledge of a controller that stabilizes every candidate model. The additional expected cost of selecting stabilizing controllers and handling incorrect selections is at most logarithmic, preserving the logarithmic and sublinear polynomial regret guarantees.

The remainder of the paper is organized as follows. Section~\ref{sec:model} formulates the problem. Section~\ref{sec:results} presents the regret classification, attaining policy and further consequences. Sections~\ref{sec:learner} and~\ref{sec:global} establish the local and global upper bounds, respectively. Section~\ref{sec:lower} proves the matching lower bounds. Section~\ref{sec:examples} provides examples of the classification, and Section~\ref{sec:conclusion} concludes the paper. The appendices provide auxiliary proofs and state-space realizations.

\emph{Notation.} $\id_d$ is the $d\times d$ identity, with dimension
omitted when clear. Vector norms are Euclidean, matrix norms are operator
norms, and $\|\cdot\|_F$ is the Frobenius norm.
For symmetric $M,N$, $M\preceq N$ or $N\succeq M$ means that $N-M$
is positive semidefinite, and $M\succ0$ means that $M$ is positive definite.
The operators $\tr,\Cov,\diag$ denote trace, covariance and block diagonal
formation, while $\vecop,\col$ stack matrix columns and vectors, respectively.
For a set $S$, $\operatorname{int}S$ and $\overline S$ denote its interior
and closure, respectively.
For random vectors, $\|X\|_{L^r}=(\E\|X\|^r)^{1/r}$, $r\ge1$.
For matrix tuples, $\|(Z_k)_k\|_{\rm mom}^2=\sum_k\|Z_k\|_F^2$.
For causal $H(\zeta)=\sum_{k\ge0}H_k\zeta^{-k}$, set
$\|H\|_{H_2}^2=\sum_{k\ge0}\|H_k\|_F^2$ and
$\|H\|_{\Wien}=\sum_{k\ge0}\|H_k\|$.
For a real stationary process with absolutely summable covariances
$\{\Gamma_k\}_{k\in\mathbb Z}$, use $S(\zeta)=\sum_{k\in\mathbb Z}\Gamma_k\zeta^{-k}$,
$|\zeta|=1$, with $\Gamma_{-k}=\Gamma_k^\top$.
Spectral Wiener norms sum over $k\in\mathbb Z$.
At each frequency, $H^*$ denotes conjugate transpose.
Write $\|z\|_S^2=z^\top Sz$ for $S\succ0$.
The logarithm $\log$ is natural, and
$\log_2^+a=\max\{0,\log_2a\}$ for $a>0$.
For nonnegative $f,g$, $f\lesssim g$ means $f\le Cg$ for a constant
$C>0$, with $\gtrsim$ reversed and $\asymp$ denoting both bounds.
We write $f\sim g$ for $f/g\to1$.
The symbols $O,\Omega,\Theta$ denote asymptotic upper, lower and two-sided
bounds, respectively. $\widetilde O$ suppresses logarithmic factors.
For probability measures, $\wedge$ denotes their common part, whose density
relative to their sum is the minimum of the two densities.
Generic $c,C>0$ denote constants in lower and upper bounds, respectively.
They may vary between estimates and are independent of $T$, with other
dependence specified locally.

\section{Problem Formulation}\label{sec:model}
\subsection{System model and quadratic cost}\label{sec:system}
We consider a family of discrete-time linear systems indexed by an unknown parameter $\theta\in\K$.
\begin{equation}\label{eq:model}
 x_{t+1}=A_\theta x_t+B_\theta u_t+w_{t+1},\qquad
 y_t=C_\theta x_t+v_t.
\end{equation}
Here $x_t\in\mathbb R^n$, $u_t\in\mathbb R^m$, and $y_t\in\mathbb R^p$
are the state, control input, and measured output, respectively.
The process noise $w_{t+1}\sim N(0,W)$ and measurement noise
$v_t\sim N(0,V)$ are mutually independent over time and between
the two sequences.

The unknown parameter is
$\theta=(\vartheta,W,V)\in\K=\Theta\times\mathcal C_W\times\mathcal C_V$.
Since the covariance matrices are symmetric, we use $\vartheta$ and the
entries on and above the diagonal of each unknown covariance matrix as
coordinates. Interiors and closures of model sets are taken in this
coordinate space.
The set $\Theta\subset\mathbb R^d$ is a known compact ball with nonempty
interior. The matrices $A_\theta=A_\vartheta$, $B_\theta=B_\vartheta$,
and $C_\theta=C_\vartheta$ are known affine functions of $\vartheta$.
Every model is stabilizable and detectable.
For known positive covariance bounds, set
\begin{equation}\label{eq:covrange}
\mathcal C_W=\{W=W^\top:\underline w\id_n\preceq W\preceq\overline w\id_n\},\qquad \mathcal C_V=\{V=V^\top:\underline v\id_p\preceq V\preceq\overline v\id_p\}.
\end{equation}
If a covariance is known, its admissible set contains only that matrix.

The unobserved initial state satisfies $x_0\sim N(\mu_\theta,\Sigma_\theta)$,
independently of the noises. Its mean and covariance are Borel measurable
functions of $\theta$ and satisfy
$\norm{\mu_\theta}+\norm{\Sigma_\theta}^{1/2}\le B_0$ uniformly over $\K$
for some finite $B_0$. The learner is not given these functions or $B_0$.
The stage cost, with known weights $Q,R$, is
\[
 c_t=x_t^\top Qx_t+u_t^\top Ru_t,\qquad Q,R\succ0.
\]

\subsection{Causal policies and regret}
At time $t$, the controller observes $y_t$ before choosing $u_t$.
The state $x_t$ and realized stage cost $c_t$ are unobserved.
The controller may use a random variable $\omega$ to randomize its inputs.
The law of $\omega$ does not depend on $\theta$, and $\omega$ is
independent of $x_0$ and the noise sequences. Define
\[
 \F_t=\sigma(\omega,y_0,\dots,y_t,u_0,\dots,u_{t-1}).
\]
A causal policy $\pi$ consists of measurable decision rules independent
of $\theta$, with $u_t$ measurable with respect to $\F_t$.
These rules may depend on the known horizon $T$.

Let $C_T=\sum_{t<T}c_t$, and let $V_T^*(\theta)$ be the minimum expected
cumulative cost over $T$ stages for a causal output feedback controller
that knows $\theta$ and the initial distribution. Writing
$\E_\theta^\pi$ for expectation under model $\theta$ and policy $\pi$,
define the expected regret and minimax regret by
\begin{equation}\label{eq:regret}
\Reg_T^\pi(\theta)=\E_\theta^\pi C_T-V_T^*(\theta),\qquad \Reg_T(\K)=\inf_\pi\sup_{\theta\in\K}\Reg_T^\pi(\theta).
\end{equation}

Our objective is to determine how $\Reg_T(\K)$ grows with $T$ and to
construct causal policies that attain this rate.

\section{Main results}\label{sec:results}
\subsection{Model comparison}\label{sec:comparison}
We characterize the growth of $\Reg_T(\K)$ by comparing differences
between optimal control laws with the information available to
distinguish models.
For a known model $\theta$, the stationary optimal output feedback
controller combines a Kalman filter and state feedback.
Let $P_\theta$ and $\Pi_\theta$ denote the stabilizing control and prediction Riccati solutions \citep[Sections~3.3 and~7.3]{andersonmoore1989}, and define
\begin{equation}\label{eq:gains}
K_\theta=(R+B_\theta^\top P_\theta B_\theta)^{-1}B_\theta^\top P_\theta A_\theta,\qquad J_\theta=\Pi_\theta C_\theta^\top(C_\theta\Pi_\theta C_\theta^\top+V)^{-1},\qquad
 L_\theta=A_\theta J_\theta.
\end{equation}
$P_\theta$ solves
\[
P_\theta=Q+A_\theta^\top P_\theta A_\theta -A_\theta^\top P_\theta B_\theta
 (R+B_\theta^\top P_\theta B_\theta)^{-1}B_\theta^\top P_\theta A_\theta,
\]
and $\Pi_\theta$ solves
\[
\Pi_\theta=W+A_\theta\Pi_\theta A_\theta^\top -A_\theta\Pi_\theta C_\theta^\top
 (C_\theta\Pi_\theta C_\theta^\top+V)^{-1}C_\theta\Pi_\theta A_\theta^\top.
\]
The stationary optimal controller is
\begin{equation}\label{eq:oracle}
\begin{split}
 \bar x_t^\theta&=q_t^\theta+J_\theta(y_t-C_\theta q_t^\theta),\quad u_t=-K_\theta\bar x_t^\theta,\\
 q_{t+1}^\theta&=(A_\theta-L_\theta C_\theta)q_t^\theta+B_\theta u_t+L_\theta y_t,
 \quad q_0^\theta=0.
\end{split}
\end{equation}
Write $\pi^\theta$ for this controller.

The plant transfer from $u$ to $y$ is
$\mathcal G_\theta(\zeta)=C_\theta(\zeta\id-A_\theta)^{-1}B_\theta$.
Let $\mathcal H_\theta$ denote the map from $y$ to $u$ defined by controller
\eqref{eq:oracle}. Their impulse response coefficients are denoted by
$\mathcal G_{\theta,k}$, $k\ge1$, and $\mathcal H_{\theta,k}$, $k\ge0$, where
$\mathcal G_{\theta,k}=C_\theta A_\theta^{k-1}B_\theta$ and
$\mathcal H_{\theta,0}=-K_\theta J_\theta$.

The difference of two plant transfers has a realization with at most
$2n$ states. By the Cayley--Hamilton theorem, equality of their first
$2n$ impulse response coefficients implies equality of the transfers. We collect these coefficients in
\begin{equation}\label{eq:raw}
 \psi(\theta)=\col(\vecop\mathcal G_{\theta,1},\ldots,\vecop\mathcal G_{\theta,2n}).
\end{equation}

For $\theta,\eta\in\K$, define
\begin{equation}\label{eq:g-b}
g_{\theta\eta}=\sum_{k=0}^{2n}
   \norm{\mathcal H_{\theta,k}-\mathcal H_{\eta,k}}_F^2,\qquad b_{\theta\eta}=\sum_{k=1}^{2n}
   \norm{\mathcal G_{\theta,k}-\mathcal G_{\eta,k}}_F^2.
\end{equation}
Thus $g_{\theta\eta}$ measures the difference between the
optimal laws, while $b_{\theta\eta}=\|\psi(\theta)-\psi(\eta)\|^2$
measures the difference between the plant input responses.

We compare stationary residual covariances under the same controller
on sets $\A\subset\K$. On each set, $\pi^\nu$ stabilizes model $\theta$
for all $\theta,\nu\in\A$.

Fix $\nu\in\A$ and apply $\pi^\nu$ to model $\alpha\in\A$.
The prediction residual is ${\widetilde y}_t=y_t-C_\nu q_t^\nu$.
Its stationary covariance at lag $k$ is
\[
 \Gamma_k^\nu(\alpha)=\lim_{t\to\infty}
 \E_\alpha^{\pi^\nu}[{\widetilde y}_t{\widetilde y}_{t-k}^\top],\qquad k\ge0.
\]
The difference between the residual covariance sequences of two closed
loops, each with at most $2n$ states, satisfies a recurrence of order
at most $4n$ at positive lags.
For $\theta,\eta,\nu\in\A$, set
\begin{equation}\label{eq:d}
 d^\nu_{\theta\eta}=\sum_{k=0}^{4n}
 \norm{\Gamma_k^\nu(\theta)-\Gamma_k^\nu(\eta)}_F^2,
 \qquad d_{\theta\eta}=d^\theta_{\theta\eta}.
\end{equation}
In $d_{\theta\eta}$, both models use $\pi^\theta$. Reversing
$\theta$ and $\eta$ changes the controller to $\pi^\eta$.
We combine the covariance and input response differences as
\begin{equation}\label{eq:info}
 \mathcal I_{\theta\eta}=d_{\theta\eta}
                       +g_{\theta\eta}b_{\theta\eta}.
\end{equation}
\begin{definition}\normalfont
We call $\A$ \emph{compatible} if any two models in $\A$ with
different optimal control laws differ in their residual covariances
or input responses. Equivalently, $\mathcal I_{\theta\eta}=0$ implies
$g_{\theta\eta}=0$ for all $\theta,\eta\in\A$.
\end{definition}

If $\A$ is compatible and $g_{\theta\eta}>0$ for some
$\theta,\eta\in\A$, define
\begin{equation}\label{eq:rho}
 \rho_\A=\min\{r\ge1:\ \exists c>0,\ \mathcal I_{\theta\eta}\ge
            c g_{\theta\eta}^r\text{ for all }\theta,\eta\in\A\}.
\end{equation}
The exponent $\rho_\A$ describes how rapidly $\mathcal I_{\theta\eta}$
can vanish as $g_{\theta\eta}$ approaches zero.

A set is \emph{semialgebraic} if it can be described by a finite
combination of polynomial equalities and inequalities.

For a set $\A\subset\K$, consider the following conditions.
\begin{enumerate}[label=(C\arabic*),ref=(C\arabic*),leftmargin=*]
\item\label{cond:cover-geometry} $\A$ is compact and semialgebraic, with
$\A=\overline{\operatorname{int}\A}$.
\item\label{cond:cover-stability} $\pi^\nu$ stabilizes model $\theta$
for all $\theta,\nu\in\A$.
\end{enumerate}
Let $\mathfrak C(\K)$ denote the collection of finite covers
$\K=\bigcup_{i=1}^M\A_i$, fixed independently of $T$, whose members
satisfy conditions~\ref{cond:cover-geometry}--\ref{cond:cover-stability}
and for which the following holds.
\begin{enumerate}[label=(C\arabic*),ref=(C\arabic*),leftmargin=*,start=3]
\item\label{cond:cover-close} For some $\delta>0$, every pair
$\theta,\eta\in\K$ with $b_{\theta\eta}\le\delta$ belongs to a common $\A_i$.
\end{enumerate}

\subsection{Regret classification and exact computation}
\begin{theorem}\label{thm:exponent}
Under the assumptions in Section~\ref{sec:system}, the following hold.
\begin{enumerate}[label=(\roman*),leftmargin=*]
\item For any compatible $\A\subset\K$ satisfying
conditions~\ref{cond:cover-geometry}--\ref{cond:cover-stability},
if $g_{\theta\eta}>0$ for some $\theta,\eta\in\A$,
the minimum in \eqref{eq:rho} is attained and $\rho_\A\in\mathbb Q$.
\item $\mathfrak C(\K)\neq\varnothing$.
Suppose $g_{\theta\eta}>0$ for some $\theta,\eta\in\K$.
For any cover $\{\A_i\}_{i=1}^M\in\mathfrak C(\K)$ whose members are
compatible, define
\[
 \rho_* =\max_{\substack{1\le i\le M\\
 \exists\,\theta,\eta\in\A_i:\ g_{\theta\eta}>0}}\rho_{\A_i}.
\]
Then $\rho_*\in\mathbb Q$ and is independent of the chosen cover.
\end{enumerate}
\end{theorem}

The cover is constructed in Appendix~\ref{app:cover}, and the exponent
properties are proved in Appendix~\ref{app:geometry}.

For $\theta,\eta\in\K$ with $b_{\theta\eta}=0$,
Condition~\ref{cond:cover-close} places both models in a common covering set.
Hence $d_{\theta\eta}$ in \eqref{eq:d} is well defined on these pairs,
independently of the cover.

\begin{theorem}\label{thm:main}
Under the assumptions in Section~\ref{sec:system}, there exists a
constant $C>0$, independent of $T$, such that the following hold.
\begin{enumerate}[label=(\roman*),leftmargin=*]
\item If $g_{\theta\eta}=0$ for all $\theta,\eta\in\K$, then
$0\le\Reg_T(\K)\le C$.
\item If there exist $\theta,\eta\in\K$ with
$b_{\theta\eta}=d_{\theta\eta}=0<g_{\theta\eta}$, then
$\Reg_T(\K)=\Theta(T)$.
\item In all other cases, $\rho_*$ in
Theorem~\ref{thm:exponent}(ii) is well defined, and
\begin{equation}\label{eq:rate}
 \Reg_T(\K)=\begin{cases}
 \Theta(\log(eT)),&\rho_*=1,\\
 \Theta(T^{1-1/\rho_*}),&\rho_*>1.
 \end{cases}
\end{equation}
\end{enumerate}
For each horizon $T\ge1$, there exists a causal policy $\pi$
attaining the corresponding regret upper bound.
\end{theorem}

A larger $\rho_*$ allows $\mathcal I_{\theta\eta}$ to vanish
faster relative to $g_{\theta\eta}$, leading to faster regret growth.
For $\rho_*>1$, the regret exponent $1-1/\rho_*$ is also rational.
For example, the family in Section~\ref{sec:cuberoot} has $\rho_*=3/2$
and $\Reg_T(\K)=\Theta(T^{1/3})$.
For model families and costs specified by real algebraic data,
the following proposition gives a finite procedure for determining the
regret class and, in case (iii), the exact $\rho_*$ before observations
are collected.

\begin{proposition}\label{prop:compute}
Suppose the coefficients of the affine maps $A_\theta,B_\theta,C_\theta$,
the entries of $Q,R$, and the parameters defining $\K$ are real algebraic
numbers given exactly. A terminating procedure determines which case
of Theorem~\ref{thm:main} holds and constructs a cover
$\{\A_i\}_{i=1}^M\in\mathfrak C(\K)$. In case (iii), it returns the exact rational
values of $\rho_{\A_i}$ for each member with varying optimal law and
their maximum $\rho_*$. The procedure requires neither the
initial distribution maps nor their bound $B_0$.
\end{proposition}

The algebraic inputs are specified by defining polynomials and isolating
intervals. Rational inputs, including exactly specified finite decimals,
satisfy this requirement.

The procedure encodes Riccati solutions by polynomial equations
and Schur Lyapunov inequalities, and moments by stationary covariance
equations. Quantifier elimination decides the
stability and separation tests \citep[Theorem~2.27]{basu14} and
the candidate rational exponents in Appendix~\ref{app:compute}.

\begin{remark}\normalfont
For
\[
 x_{t+1}=u_t+w_{t+1},\qquad y_t=x_t+v_t,
\]
with $Q=R=1$ and unknown $W,V\in[1,2]$, we have $K_\theta=0$
for every model. The optimal input is therefore $u_t=0$,
regardless of the unknown covariances, and case (i) applies.

For case (ii), consider a two-dimensional system with state
$x_t=(z_t,h_t)^\top$ and process noise $w_t=(w_t^z,w_t^h)^\top$.
Only the first state component is measured, with additive noise.
The system is
\[
z_{t+1}=w^z_{t+1},\qquad y_t=z_t+v_t,\qquad h_{t+1}=\theta z_t+u_t+w^h_{t+1},
\]
with unit noise covariances and identity cost weights. In stationarity,
$\E_\theta[z_t\mid\F_t]=y_t/2$, and minimizing the input cost together
with the next-state cost gives $u_t=-\theta y_t/4$.
Since the input does not affect $y_t$ and the prediction residual is $y_t$,
we have $b_{\theta\eta}=d_{\theta\eta}=0$. Thus, for $\theta\ne\eta$,
\[
 \mathcal I_{\theta\eta}=0,\qquad
 g_{\theta\eta}=\frac{(\theta-\eta)^2}{16}>0,
\]
and case (ii) applies.
\end{remark}

In the preceding example, the unobserved state $h_t$ enters the cost.
When the cost depends only on the output and input, linear regret is
excluded and the remaining rates satisfy the following restriction.

\begin{corollary}\label{cor:output}
Replace the stage cost by $y_t^\top Q_y y_t+u_t^\top Ru_t$, where
$Q_y\succ0$ and $R\succ0$ are known and do not depend on $\theta$.
If the stationary optimal law is constant on $\K$, then $\Reg_T(\K)$
is bounded. Otherwise, $1\le\rho_*\le2$, and
\[
\Reg_T(\K)=
\begin{cases}
\Theta(\log(eT)),&\rho_*=1,\\
\Theta(T^\alpha),\quad\alpha=1-1/\rho_*,&1<\rho_*\le2.
\end{cases}
\]
Here $\rho_*$ is defined as in Theorem~\ref{thm:exponent} using the
optimal control laws for the output cost, and $\alpha$ is rational
with $0<\alpha\le\frac12$ in the second case.
\end{corollary}

\begin{remark}
The output cost applies the same known weights to $y_t$ and $u_t$
under every model. Its quadratic structure and optimality give,
on each covering set,
\[
 g_{\theta\eta}\lesssim d_{\theta\eta}+b_{\theta\eta},\qquad
 \mathcal I_{\theta\eta}\gtrsim g_{\theta\eta}^{2}.
\]
The second inequality gives $\rho_\A\le2$ on each set with varying
optimal law. Hence $\rho_*\le2$ and $\alpha=1-1/\rho_*\le1/2$
in the power case. Appendix~\ref{app:output} proves the comparison.
\end{remark}

\subsection{Control policy}\label{sec:attaining}
The policy attains the upper bounds in Theorem~\ref{thm:main}.
When the optimal law varies over $\K$, use pulse observations to estimate
$\psi(\theta)$ and select a covering set $\A$. If an output or proposed input exceeds its norm
threshold, reselect $\A$ and restart (Section~\ref{sec:global}).

On the selected compatible set $\A$ with varying optimal laws,
Algorithm~\ref{alg:local} uses $\rho_\A$ from \eqref{eq:rho} and
updates $\nu_j\in\A$ between blocks. Block $j$ uses random input $\sigma_j e_t$.
Choose $N_0$ to reduce the initial state's contribution by at least
half per full block (Lemma~\ref{lem:block}). Block lengths
$N_j$ and amplitudes $\sigma_j\ge0$ satisfy
\begin{equation}\label{eq:schedule}
 N_j=N_0 2^j,\qquad \sigma_j^2=(N_0/N_j)^{1/\rho_\A},\quad j\ge0.
\end{equation}
Use independent vectors $e_t\in\{-1,1\}^m$ with independent
Rademacher coordinates, independent of $x_0$ and the noise sequences $\{w_t\}_{t\ge1}$ and $\{v_t\}_{t\ge0}$.
The control recursion in block $j$ is
\begin{equation}\label{eq:implement}
\begin{aligned}
 {\widetilde y}_t&=y_t-C_{\nu_j} q_t, &\bar x_t^{\nu_j}&=q_t+J_{\nu_j}{\widetilde y}_t,\\
 u_t&=-K_{\nu_j}\bar x_t^{\nu_j}+\sigma_j\,e_t,
 &q_{t+1}&=A_{\nu_j}\bar x_t^{\nu_j}+B_{\nu_j} u_t.
\end{aligned}
\end{equation}
\begin{localpolicy}\label{alg:local}\leavevmode
\begin{enumerate}[label=\arabic*),leftmargin=*,nosep]
\item At the start $b_0$ of the local run, choose $\nu_0\in\A$
and set $j=0$.
\item Set $N_j$ and $\sigma_j$ by \eqref{eq:schedule},
$b_{j+1}=b_j+N_j$, and reset $q_{b_j}=0$.
\item For $b_j\le t<\min(b_{j+1},T)$, apply \eqref{eq:implement}.
End the local run at horizon $T$ or when the output or proposed input exceeds
its norm threshold (Section~\ref{sec:global}).
\item If the block completes without stopping, fit $\nu_{j+1}$ by
least squares \eqref{eq:project_fit} to the residual autocovariances and
cross covariances with probes \eqref{eq:sample}.
Increase $j$ by one and return to step~2.
\end{enumerate}
\end{localpolicy}
Each block incurs $O(N_j\sigma_j^2+1)$ regret
(Proposition~\ref{prop:local}). For $\rho_\A=1$, \eqref{eq:schedule} gives
$N_j\sigma_j^2=N_0$, so $O(\log(eT))$ blocks give $O(\log(eT))$ regret.
For $\rho_\A>1$, geometric growth of $N_j\sigma_j^2$ gives $O(T^{1-1/\rho_\A})$.

A common optimal controller on $\K$ gives $O(1)$ regret
(Lemma~\ref{lem:block}). Otherwise, use the selected set's common
controller if its optimal law is constant. For varying laws, use
Algorithm~\ref{alg:local} if compatible and stabilizing feedback
with $O(T)$ regret otherwise.
Repeated localization adds $O(\log(eT))$ to the worst local bound
(Proposition~\ref{prop:global}), determined by $\rho_*$ in case (iii).
The matching lower bounds in Section~\ref{sec:lower} then yield \eqref{eq:rate}.

\section{Model estimation and local regret}\label{sec:learner}
The update fits theoretical moments to sample estimates.

\subsection{Prediction residuals}
For $\theta,\nu\in\A$ and $\sigma\in[0,1]$, set
$F_{c,\nu}=A_\nu-B_\nu K_\nu$,
$F_{o,\nu}=A_\nu-L_\nu C_\nu$ and $F_\nu=\id-J_\nu C_\nu$.
For $Z_t=\col(x_t,q_t)$, the closed loop is
\begin{equation}\label{eq:physical_real}
Z_{t+1}=M_{\theta,\nu}Z_t+E_ww_{t+1} +E_{v,\theta,\nu}v_t+\sigma E_{e,\theta,\nu}e_t,\qquad {\widetilde y}_t=H_{\theta,\nu}Z_t+v_t,
\end{equation}
where
\begin{equation}\label{eq:phys_matrix}
\begin{gathered}
M_{\theta,\nu}=\begin{pmatrix}
 A_\vartheta-B_\vartheta K_\nu J_\nu C_\vartheta&-B_\vartheta K_\nu F_\nu\\
 F_{c,\nu}J_\nu C_\vartheta&F_{c,\nu}F_\nu
 \end{pmatrix},\\
E_w=\binom{\id_n}{0},\quad
 E_{v,\theta,\nu}=\binom{-B_\vartheta K_\nu J_\nu}{F_{c,\nu}J_\nu},\qquad E_{e,\theta,\nu}=\binom{B_\vartheta}{B_\nu},\quad
 H_{\theta,\nu}=(C_\vartheta,-C_\nu).
\end{gathered}
\end{equation}
Compactness and condition~\ref{cond:cover-stability} make $M_{\theta,\nu}$ uniformly Schur.
Write $\theta=(\vartheta,\Xi)$ with $\Xi=\diag(W,V)$, and set
$\mathcal C=\{\diag(W,V):W\in\mathcal C_W,V\in\mathcal C_V\}$ and
$E_{\theta,\nu}=(E_w,E_{v,\theta,\nu})$.
The stationary covariances $S_{\nu,\sigma^2}(\theta)=\Cov(Z_t)$ and
$G_{\nu,\sigma^2}(\theta)=\Cov(Z_{t+1},{\widetilde y}_t)$ satisfy
\begin{equation}\label{eq:stein}
\begin{gathered}
S_{\nu,\sigma^2}(\theta)
 =M_{\theta,\nu}S_{\nu,\sigma^2}(\theta)M_{\theta,\nu}^\top +E_{\theta,\nu}\Xi E_{\theta,\nu}^\top
       +\sigma^2 E_{e,\theta,\nu}E_{e,\theta,\nu}^\top,\\
G_{\nu,\sigma^2}(\theta)
 =M_{\theta,\nu}S_{\nu,\sigma^2}(\theta)H_{\theta,\nu}^\top+E_{v,\theta,\nu}V.
\end{gathered}
\end{equation}
The cross term is $E_{v,\theta,\nu}V=\Cov(Z_{t+1},v_t)$.
Collect the residual autocovariances at lags $0,\ldots,4n$ and cross
covariances with probes at lags $1,\ldots,2n$ as
\begin{equation}\label{eq:feature}
\Phi_{\nu,\sigma^2}(\vartheta,\Xi)=
 \bigl(H_{\theta,\nu}S_{\nu,\sigma^2}(\theta)H_{\theta,\nu}^\top+V, (H_{\theta,\nu}M_{\theta,\nu}^{k-1}G_{\nu,\sigma^2}(\theta))_{k=1}^{4n}, (\sigma H_{\theta,\nu}M_{\theta,\nu}^{k-1}E_{e,\theta,\nu})_{k=1}^{2n}\bigr).
\end{equation}
For fixed dynamics, both $S_{\nu,\sigma^2}$ and $\Phi_{\nu,\sigma^2}$
are affine in $\Xi$.

\begin{lemma}\label{lem:feature}
Uniformly for $\theta,\eta,\nu\in\A$ and $\sigma\in[0,1]$,
\begin{equation}\label{eq:feature_equiv}
 \|\Phi_{\nu,\sigma^2}(\theta)-\Phi_{\nu,\sigma^2}(\eta)\|_{\rm mom}^2
 \asymp d^\nu_{\theta\eta}+\sigma^2 b_{\theta\eta}.
\end{equation}
\end{lemma}
\begin{proof}
Allow $r_{\nu,t}=u_t+K_\nu\bar x_t^\nu$. Write the stable inverse of
\eqref{eq:implement} as $u=\mathcal C_\nu{\widetilde y}+\mathcal D_\nu r_\nu$.
For $\alpha\in\A$, let $\mathcal T_{\alpha,\nu}$ be the transfer
$r_\nu\mapsto{\widetilde y}$ in \eqref{eq:physical_real}, with $\sigma e=r_\nu$.
With $\mathcal R_\nu=\id-C_\nu(\zeta\id-F_{o,\nu})^{-1}L_\nu$,
the inverse and plant equations give
\begin{equation}\label{eq:short_probe_identity}
\begin{gathered}
\mathcal C_\nu=-\zeta K_\nu(\zeta\id-F_{c,\nu})^{-1}J_\nu,\qquad \mathcal D_\nu=\id-K_\nu(\zeta\id-F_{c,\nu})^{-1}B_\nu,\\
\mathcal T_{\alpha,\nu}=\mathcal R_\nu(\mathcal G_\alpha-\mathcal G_\nu)
 (\mathcal D_\nu+\mathcal C_\nu\mathcal T_{\alpha,\nu}).
\end{gathered}
\end{equation}
Since $\mathcal R_\nu,\mathcal D_\nu$ have identity feedthrough and
$\mathcal G_\alpha-\mathcal G_\nu$ is strictly causal, the
first $2n$ coefficient map is triangular and uniformly invertible:

\begin{equation}\label{eq:short_b}
 \sum_{k=1}^{2n}\|\mathcal T_{\theta,\nu,k}-
                  \mathcal T_{\eta,\nu,k}\|_F^2
 \asymp b_{\theta\eta}.
\end{equation}
Let $A=(\Gamma_k^\nu(\theta)-\Gamma_k^\nu(\eta))_{k=0}^{4n}$.
Let $B$ be the covariance components of
$\Phi_{\nu,1}(\theta)-\Phi_{\nu,1}(\eta)$ minus $A$, and set
$T=(\mathcal T_{\theta,\nu,k}-\mathcal T_{\eta,\nu,k})_{k=1}^{2n}$.
For $\Delta\Phi=\Phi_{\nu,\sigma^2}(\theta)-\Phi_{\nu,\sigma^2}(\eta)$,
the Stein equation gives
$\Delta\Phi=(A+\sigma^2B,\sigma T)$ and
$\|A\|_{\rm mom}^2=d^\nu_{\theta\eta}$.
Subtracting the probe spectra
$\mathcal T_{i,\nu}\mathcal T_{i,\nu}^*$, $i=\theta,\eta$, and using
\eqref{eq:finite_transfer} gives $\|B\|_{\rm mom}\le C\|T\|_{\rm mom}$.
Since $\sigma^4\le\sigma^2$, the triangle inequality in both directions yields
\[
\|\Delta\Phi\|_{\rm mom}^2
 =\|A+\sigma^2B\|_{\rm mom}^2+\sigma^2\|T\|_{\rm mom}^2 \asymp\|A\|_{\rm mom}^2+\sigma^2\|T\|_{\rm mom}^2
 \asymp d^\nu_{\theta\eta}+\sigma^2 b_{\theta\eta}.
\]
\end{proof}

\subsection{Covariance estimates and model fitting}\label{sec:project_fit}
For a completed block of length $N>8n+2$ beginning at $b$,
use ${\widetilde y}_t=e_t=0$ for $t<b$ in the statistics
\begin{equation}\label{eq:sample}
\widehat\Gamma_k=N^{-1}\sum_{t=b}^{b+N-1}{\widetilde y}_t{\widetilde y}_{t-k}^\top,
  0\le k\le4n,\qquad \widehat H_k=N^{-1}\sum_{t=b}^{b+N-1}{\widetilde y}_t e_{t-k}^\top,
  1\le k\le2n.
\end{equation}
Set $\widehat\Phi=((\widehat\Gamma_k)_{k=0}^{4n},
(\widehat H_k)_{k=1}^{2n})$, also denoted $\widehat\Phi_b$ when
marking its block start.

Use the response sets $\A$ of Appendix~\ref{app:cover}, refined for
$\rho_\A=1$. They satisfy $\A=\Theta_\A\times\mathcal C$, where
$\Theta_\A$ is the dynamics projection, since $\psi$ depends only on $\vartheta$.
Discretize $\Theta_\A$ by closed dyadic cubes of diameter $\le N^{-1/2}$.
From each intersecting cube, take the lexicographically first point
of $\Theta_\A$. Order them lexicographically as $\mathcal N_N$ and fit by
\begin{equation}\label{eq:project_fit}
(\widehat\vartheta,\widehat\Xi)
 =\operatorname*{argmin}_{\vartheta\in\mathcal N_N,
                                      \Xi\in\mathcal C}
 \|\widehat\Phi-\Phi_{\nu,\sigma^2}(\vartheta,\Xi)\|_{\rm mom}^2,\qquad \widehat\theta=(\widehat\vartheta,\widehat W,\widehat V).
\end{equation}
Here $\widehat\Xi=\diag(\widehat W,\widehat V)$.
Set
\[
 \Xi_c=\diag\!\left(\tfrac{\underline w+\overline w}{2}\id_n,
                         \tfrac{\underline v+\overline v}{2}\id_p\right),
\]
using the known covariance for any singleton block.
Among minimizers, choose the covariance nearest $\Xi_c$ in Frobenius
norm, then the first tied grid point. Each covariance solution set is
compact and convex, so this choice is unique. Continuity ensures Borel selection.

\begin{lemma}\label{lem:profile}
The update \eqref{eq:project_fit} is measurable and, uniformly
for $\theta,\nu\in\A$, obeys
\begin{equation}\label{eq:fit_stability}
\|\Phi_{\nu,\sigma^2}(\widehat\theta)-\Phi_{\nu,\sigma^2}(\theta)\|_{\rm mom}^2 \le8\|\widehat\Phi-\Phi_{\nu,\sigma^2}(\theta)\|_{\rm mom}^2+C/N.
\end{equation}
\end{lemma}
\begin{proof}
Choose a grid point $\vartheta_N$ within $N^{-1/2}$ of
$\vartheta$, retaining the true $\Xi$. Lipschitz continuity gives
$\|\Phi_{\nu,\sigma^2}(\vartheta_N,\Xi)-\Phi_{\nu,\sigma^2}(\theta)\|_{\rm mom}\le L/\sqrt N$.
Set $\Delta=\|\widehat\Phi-\Phi_{\nu,\sigma^2}(\theta)\|_{\rm mom}$.
Optimality and two triangle inequalities yield
\[
\|\Phi_{\nu,\sigma^2}(\widehat\theta)-\Phi_{\nu,\sigma^2}(\theta)\|_{\rm mom} \le2\Delta+L/\sqrt N.
\]
Squaring with $(a+b)^2\le2(a^2+b^2)$ proves
\eqref{eq:fit_stability}.
\end{proof}

\begin{lemma}\label{lem:block}
Compactness and condition~\ref{cond:cover-stability} give $\|M_{\theta,\nu}^k\|
\le C_0\gamma^k$, $\gamma<1$. Choose $N_0>8n+2$ with
$C_0\gamma^{N_0}\le\tfrac12$. Then
\begin{equation}\label{eq:reset_contraction}
 \sup_{\theta,\nu\in\A}\sup_{N\ge N_0}
 \|(\id_n,0)M_{\theta,\nu}^{N}(\id_n,0)^\top\|\le\tfrac12.
\end{equation}
Assume uniformly bounded moments of each fixed order at the start
of the local run. With $q_b=0$ and controllers selected from the
information at block starts, these bounds persist within and across
blocks of length $\ge N_0$. For every integer $q\ge1$,
\begin{equation}\label{eq:sample_moment}
 \E_\theta\|\widehat\Phi_b-\Phi_{\nu,\sigma^2}(\theta)\|_{\rm mom}^{2q}
 \le C_qN^{-q}.
\end{equation}
With $J_\theta^*=\lim_{N\to\infty}N^{-1}\E_\theta^{\pi^\theta}C_N$
the stationary optimal cost per stage,
\begin{equation}\label{eq:blockcost}
 \E_\theta[\text{block cost}]-NJ_\theta^*
 \le CN\{\E_\theta g_{\theta\nu}+\sigma^2\}+C.
\end{equation}
For arbitrary block initialization, the cost bound has an additional
$C(1+\E_\theta\|x_b\|^2)$ term and holds for truncated blocks. Uniformly,
$V_T^*(\theta)=TJ_\theta^*+O(1)$.
\end{lemma}
\begin{proof}
Equations \eqref{eq:blockrec} and \eqref{eq:quadratic_moment}
yield \eqref{eq:sample_moment}. Square completion, initialization
and Riccati convergence in Appendix~\ref{app:moment} yield \eqref{eq:blockcost}.
\end{proof}

\subsection{Local regret bound}

\begin{proposition}\label{prop:local}
For such a compatible $\A$ with varying optimal laws and the initial
moment bounds of Lemma~\ref{lem:block}, the learner
\eqref{eq:implement}, \eqref{eq:project_fit}, \eqref{eq:schedule}
satisfies uniformly
\begin{equation}\label{eq:localrate}
 \sup_{\theta\in\A}\Reg_T^\pi(\theta)\le
 \begin{cases}C\log(eT),&\rho_\A=1,\\CT^{1-1/\rho_\A},&\rho_\A>1.\end{cases}
\end{equation}
\end{proposition}
\begin{proof}
For completed block $j$, write $N=N_j$, $\sigma^2=\sigma_j^2$, $\nu=\nu_j$
and $\nu^+=\nu_{j+1}$. The current and updated control-law errors are
$H=g_{\theta\nu}$ and $G=g_{\theta\nu^+}$. Set the statistical error term
\[
 U=C\bigl(\|\widehat\Phi-\Phi_{\nu,\sigma^2}(\theta)\|_{\rm mom}^2+N^{-1}\bigr).
\]
Equations \eqref{eq:feature_equiv}, \eqref{eq:fit_stability} and \eqref{eq:sample_moment} give
$d^\nu_{\theta\nu^+}+\sigma^2 b_{\theta\nu^+}\le U$ and
$\E U^q\le C_qN^{-q}$ for all fixed $q>0$.
For $\rho_\A>1$, \eqref{eq:rho} and Lemma~\ref{lem:transport} give
\[
 G^{\rho_\A}\le C\{d^\nu_{\theta\nu^+}+(H+G)b_{\theta\nu^+}\}
 \le CU(1+(H+G)/\sigma^2).
\]
If $G>(H+\sigma^2)/8$, the last factor is at most $9G/\sigma^2$.
Otherwise $G\le(H+\sigma^2)/8$. In either case,
\[
 G\le(H+\sigma^2)/8+(9CU/\sigma^2)^{1/(\rho_\A-1)}.
\]
Since $N\sigma^{2\rho_\A}=N_0$, the last term has expectation at most
$C(N\sigma^2)^{-1/(\rho_\A-1)}\le C\sigma^2$.
For $\rho_\A=1$, the refinement in Appendix~\ref{app:endpoint-refinement} gives
\begin{equation}\label{eq:endpoint}
 g_{\theta\eta}\le\tfrac18g_{\theta\nu}+Cd^\nu_{\theta\eta}
 \qquad\theta,\eta,\nu\in\A.
\end{equation}
Thus $G\le H/8+CU$ and $\E U\le C/N\le C\sigma^2$.
Both cases, with $\bar g_j=\E g_{\theta\nu_j}$, yield
\begin{equation}\label{eq:contract}
 \bar g_{j+1}\le\tfrac18\bar g_j+C\sigma_j^2,\qquad \bar g_j\le C\sigma_j^2.
\end{equation}
The second bound follows by induction from bounded $\bar g_0$ and
$1/8<\sigma_{j+1}^2/\sigma_j^2=2^{-1/\rho_\A}$. Lemma~\ref{lem:block}, including
the truncated block, gives for $T\ge N_0$
\[
\Reg_T^\pi(\theta)\le C\sum_{j\le\log_2(eT/N_0)}(N_j\sigma_j^2+1).
\]
Here $N_j\sigma_j^2=N_0$ for $\rho_\A=1$. Otherwise the sum is
$O(T^{1-1/\rho_\A})$. Smaller horizons are absorbed into $C$,
proving \eqref{eq:localrate}.
\end{proof}

\section{Global stabilization by repeated localization}\label{sec:global}
Set $R_t=\max_{k<t}(\norm{y_k}+\norm{u_k})$, with $R_0=0$.
Apply $u_t=0$ for $t<\min(n,T)$. If $T>n$, set
\begin{equation}\label{eq:threshold_init}
 a_T=M_0\sqrt{\log(eT)},\qquad h_0=a_T+R_n,
\end{equation}
where $M_0\ge1$ depends only on known Gaussian noise bounds.
Starting from $\tau_0=n$, trial $j$ uses
\begin{equation}\label{eq:trial_scales}
 h_j=2^jh_0,\quad \overline R_j=R_{\tau_j},\quad
 \Lambda_j=\overline R_j+h_j.
\end{equation}
Every trial applies $m$ coordinate pulses with $\ell=2n$ and fixed
$C_{\rm p},a_{\rm p}>1$.
The vector $e_r^{(m)}$ is the $r$th standard basis vector in $\mathbb R^m$.
\begin{equation}\label{eq:pulses}
 u_{\tau_j+r\ell}=C_{\rm p}\Lambda_ja_{\rm p}^{\,r} e^{(m)}_{r+1},
 \qquad r=0,\ldots,m-1.
\end{equation}
Other inputs are zero through $\tau_j+m\ell$, stopping at $T$.
Normalize $y_{\tau_j+r\ell+k}$, $1\le k\le\ell$, by pulse $r$'s
amplitude to estimate response column $r+1$ in \eqref{eq:raw}.
For the final cover in Appendix~\ref{app:cover}, set
$\varepsilon_0=\min_i\varepsilon_i>0$, where $4\varepsilon_i$
is the outer radius of the set centered at $a_i$.
After collection, select the first index minimizing
$\norm{\widehat\psi-\psi(a_i)}/\varepsilon_i$,
where $\widehat\psi$ is the response estimate.
At $b_0=\tau_j+m\ell+1$, reset $q=0$ and run
Section~\ref{sec:attaining} on the selected set.

During continuation, if the output or proposed input has norm greater than
$H_j=C_{\rm out}\Lambda_j$ at time $t$, use $u_t=0$, discard $q$,
and start the next trial at $\tau_{j+1}=t+1$.

For analysis, set
\begin{equation}\label{eq:noise_env}
 D_N=\max\Big\{\sqrt m,\max_{1\le t\le N}\norm{w_t},
 \max_{0\le t\le N}\norm{v_t}\Big\}.
\end{equation}
Choose known $C_{\rm p},a_{\rm p},C_{\rm out}$ from
Appendix~\ref{app:cover}'s uniform bounds so that, when collection
ends before $N$ and $D_N\le h_j$, the response error is $<\varepsilon_0$
and continuation norms are $<C_{\rm out}\Lambda_j$.
The selected $a_i$ satisfies $\norm{\psi(\theta)-\psi(a_i)}<3\varepsilon_i$,
so $\theta\in\A_i$.
For each fixed $r\ge1$ and large enough known $M_0$,
\begin{equation}\label{eq:noise_tail}
 \Prob(D_T>a_T)\le CT^{-12},\qquad
 \norm{D_T}_{L^r}\le C_r a_T.
\end{equation}
Thus every trial rejected before $N$ has $h_j<D_N$, and pathwise
\begin{equation}\label{eq:pathbound}
 \max_{t<N}(\norm{y_t}+\norm{u_t})
 \le Ch_0(1+D_N/h_0)^q
\end{equation}
for fixed integer $q$.

For a covering set $\A$, put
$G_\A=\max_{\theta,\eta\in\A}g_{\theta\eta}$.

\begin{proposition}\label{prop:global}
On compatible $\A_i$ with varying optimal laws, use
Algorithm~\ref{alg:local} and let $r_i(T)$ be the rate in Proposition~\ref{prop:local}.
On sets with a common optimal controller, use it and set $r_i(T)=1$.
On incompatible sets, fix $\pi^\nu$, $\nu\in\A_i$, and set $r_i(T)=T$.
The restart policy $\pi$ satisfies, for $T\ge1$,
\begin{equation}\label{eq:globalreg}
 \sup_{\theta\in\K}\Reg_T^\pi(\theta)
 \le C\{\log(eT)+\max_i r_i(T)\}.
\end{equation}
It satisfies
\[
\sup_{\theta\in\K,\,0\le t<T}
 \E_\theta^\pi(\norm{x_t}^2+\norm{u_t}^2)
 \le C\log(eT),\qquad \sup_{\theta\in\K}\E_\theta^\pi C_T^2
 \le CT^2\log^2(eT).
\]
If $g_{\theta\eta}=0$ for all $\theta,\eta\in\K$, the common
stationary controller $\pi$, used from time zero, satisfies these
moment bounds and $\sup_{\theta\in\K}\Reg_T^\pi(\theta)\le C$.
The constant $C>0$ is independent of $T$.
\end{proposition}
\begin{proof}
Absorb $T\le n+m\ell+1$ into $C$.
The prefix has expected cost $O(a_T^2)$ and $\|x_{b_0}\|_{L^r}=O(a_T)$.
For each $\A_i$ containing $\theta$, continue with the same prefix,
noises and seed. This agrees with the actual path on $\{D_T\le a_T\}$
when $\A_i$ is selected.
Weight the bounds by the selection indicators, determined by the
pulse record before $b_0$, and sum over $i$.

For compatible $\A=\A_i$ with varying laws, index local blocks by $j$.
The reset recursion \eqref{eq:blockrec} gives
$\|x_{b_j}\|_{L^r}\le C_r(1+2^{-j}a_T)$.
With $j_*=\lceil\log_2 a_T\rceil$, the first $j_*$ blocks have excess cost
at most $C\sum_{j<j_*}(N_j+1+4^{-j}a_T^2)=O(a_T+a_T^2)$.
After $j_*$, boundary moments are bounded.
With $\bar g_j=\E_\theta g_{\theta\nu_j}$, \eqref{eq:contract} gives
\[
 \bar g_{j_*+k}\le8^{-k}G_\A+C\sigma_{j_*+k}^2,\qquad k\ge0.
\]
Since $N_{j_*}=O(a_T)$, the initial controller error costs at most
$CN_{j_*}G_\A\sum_{k\ge0}(2/8)^k=O(a_T)$.
The remaining costs, including resets and truncation, sum to the
local bound (Proposition~\ref{prop:local}).
The prescribed common optimal and stabilizing controllers give excess
cost $O(a_T^2+r_i(T))$ by Lemma~\ref{lem:block}.
Thus each continuation has excess cost $O(\log(eT)+r_i(T))$ and,
by stability, $\sup_{t<T}(\|x_t\|_{L^r}+\|u_t\|_{L^r})=O(a_T)$.

On $\{D_T>a_T\}$, \eqref{eq:pathbound}, \eqref{eq:observerbound},
\eqref{eq:noise_tail} and Cauchy--Schwarz give $O(1)$ contributions
to expected cost and second and fourth moments of the state and input.
Combine the continuation and exceptional contributions with the
benchmark comparison in Lemma~\ref{lem:block} to obtain
\eqref{eq:globalreg} and both moment bounds,
using $C_T^2\le T\sum_{t<T}c_t^2$.
A globally common optimal controller agrees with each $\pi^\theta$,
so Lemma~\ref{lem:block} gives bounded regret and second and fourth moments.
\end{proof}

\section{Lower bounds by a two-model comparison}\label{sec:lower}
For a covering set $\A$ and $0<s\le G_\A$, define
\begin{equation}\label{eq:profile}
 D_\A(s)=\min_{\substack{\theta,\eta\in\A\\g_{\theta\eta}\ge s}}
  \{d_{\theta\eta}+s b_{\theta\eta}\}.
\end{equation}
Small $D_\A(s)$ identifies model pairs that are hard to distinguish.
Whenever $\rho_\A$ is defined, $D_\A(s)\asymp s^{\rho_\A}$ as $s\downarrow0$
(Appendix~\ref{app:geometry}).

\subsection{Regret and comparison of residuals}
Let $\widehat x_t^\theta=\E_\theta[x_t\mid\F_t]$. Starting from
$P_{\theta,0}=0$, define for integers $h\ge1$
$S_{\theta,h}=R+B_\theta^\top P_{\theta,h-1}B_\theta$,
$K_{\theta,h}=S_{\theta,h}^{-1}B_\theta^\top P_{\theta,h-1}A_\theta$
and $P_{\theta,h}=Q+A_\theta^\top P_{\theta,h-1}
(A_\theta-B_\theta K_{\theta,h})$.
Conditional completion of squares gives
\begin{equation}\label{eq:bellman}
\Reg_T^\pi(\theta)=\sum_{t<T}\E_\theta^
\pi\norm{u_t+K_{\theta,T-t}\widehat x_t^\theta}_{S_{\theta,T-t}}^2.
\end{equation}
Evaluate the predictors in \eqref{eq:oracle} on the same record $(y,u)$,
with $q_0^i=0$. Their prediction residuals and deviations from the
stationary optimal feedback are
\begin{equation}\label{eq:direct-innov}
 {\widetilde y}_{i,t}=y_t-C_iq_t^i,\qquad
 r_{i,t}=u_t+K_i(q_t^i+J_i{\widetilde y}_{i,t}).
\end{equation}
Comparison with the finite-horizon control law in
Appendix~\ref{app:energy} gives
\begin{equation}\label{eq:energy}
\sum_{t<T/2}\E_\theta^\pi\norm{r_{\theta,t}}^2
\le C(\Reg_T^\pi(\theta)+1).
\end{equation}
For $\theta,\eta\in\A$, the stable causal map between residual sequences is
\begin{equation}\label{eq:twoport}
\binom{{\widetilde y}_\eta}{r_\eta}=
\begin{pmatrix}\mathsf A_{\eta\theta}&\mathsf B_{\eta\theta}\\
\mathsf C_{\eta\theta}&\mathsf D_{\eta\theta}\end{pmatrix}
\binom{{\widetilde y}_\theta}{r_\theta}.
\end{equation}
Appendix~\ref{app:KL} gives a $2n$-state realization.
With $j\ge0$ indexing impulse response coefficients, uniformly
for some fixed $0<\alpha<1$,
\begin{equation}\label{eq:twoport-bounds}
\norm{\mathsf B_{\eta\theta}}_{\Wien}^2\le Cb_{\theta\eta},\quad
 \norm{\mathsf D_{\eta\theta,j}}\le C\alpha^j,\qquad \sum_{j=0}^{2n}\norm{\mathsf C_{\eta\theta,j}}_F^2
 \asymp\norm{\mathsf C_{\eta\theta}}_{\Wien}^2\asymp g_{\theta\eta}.
\end{equation}
Both $\mathsf A_{\eta\theta},\mathsf D_{\eta\theta}$ have identity
feedthrough. $\mathsf B_{\eta\theta}$ is strictly causal. Under $\theta$,
$\widetilde y_\theta$ is Gaussian, policy-independent, and its finite
covariance matrices are uniformly positive definite.

\subsection{Observed information and residual separation}
Let $\Prob_{i,N}^\pi$ be the joint law of
$(\omega,y_0,\ldots,y_{N-1})$ under model $i$ and policy $\pi$.
Write $\mathrm{KL}$ for Kullback--Leibler divergence, with
$\mathrm{KL}(P\Vert Q)=\int\log(dP/dQ)\,dP$ when $P$ is
absolutely continuous with respect to $Q$, and $+\infty$ otherwise.
For a fixed pair $\theta,\eta$, write
$D_{\mathrm{KL},N}(\theta,\eta)=\mathrm{KL}
(\Prob_{\theta,N}^{\pi^\theta}\Vert\Prob_{\eta,N}^{\pi^\theta})$.
\begin{lemma}\label{lem:KL}
For arbitrary causal randomized policies and $N\le T/2$,
\begin{equation*}
\mathrm{KL}(\Prob_{\theta,N}^\pi\Vert\Prob_{\eta,N}^\pi)
\le C\{Nd_{\theta\eta}+b_{\theta\eta}\Reg_T^\pi(\theta)+1\}.
\end{equation*}
\end{lemma}
\begin{proof}
Appendix~\ref{app:KL} gives
$D_{\mathrm{KL},N}(\theta,\eta)\le C(Nd_{\theta\eta}+1)$ and
\[
\mathrm{KL}\le2D_{\mathrm{KL},N}(\theta,\eta)+
Cb_{\theta\eta}\E_\theta\sum_{t<N}\|r_{\theta,t}\|^2.
\]
Apply \eqref{eq:energy}.
\end{proof}

\begin{lemma}\label{lem:common}
For fixed $K<\infty$ there exist $c_K>0$ and $0<\alpha<1$
such that, if $\mathrm{KL}(\Prob_{\theta,t+1}^\pi\Vert
\Prob_{\eta,t+1}^\pi)\le K$, then for $t\ge2n$,
\begin{equation}\label{eq:common}
c_Kg_{\theta\eta}\le C\sum_{k=0}^t\alpha^{t-k}
\{\E_\theta^\pi\norm{r_{\theta,k}}^2+
\E_\eta^\pi\norm{r_{\eta,k}}^2\}.
\end{equation}
\end{lemma}
\begin{proof}
For the same observed record,
\begin{equation}\label{eq:direct-common-identity}
X_t=(\mathsf C_{\eta\theta}{\widetilde y}_\theta)_t
=r_{\eta,t}-(\mathsf D_{\eta\theta}r_\theta)_t.
\end{equation}
Under $\theta$, $X_t$ is policy-independent Gaussian with
$\tr\Cov(X_t)\ge cg_{\theta\eta}$ by \eqref{eq:twoport-bounds}.
For $g_{\theta\eta}>0$, some coordinate has variance
$\ge cg_{\theta\eta}$, so its Gaussian density gives
$\Prob_\theta^\pi(\norm{X_t}^2\le\varepsilon g_{\theta\eta})
\le C\sqrt\varepsilon$.
The two history laws share mass $\ge e^{-K}/2$.
Choose $\varepsilon>0$ with $C\sqrt\varepsilon\le e^{-K}/4$. Then
$\int\norm{X_t}^2\,d(\Prob_{\theta,t+1}^\pi\wedge\Prob_{\eta,t+1}^\pi)
\ge c_Kg_{\theta\eta}$.
Integrating \eqref{eq:direct-common-identity} under the common measure,
bounding residual energies by their own model laws, and using
geometric decay of $\mathsf D_{\eta\theta}$ gives \eqref{eq:common}.
\end{proof}

\subsection{Regret lower bounds}
\begin{proposition}\label{prop:powerlower}
If $TD_\A(s)\le1$, then
\begin{equation}\label{eq:twopoint}
\Reg_T(\A)\ge cTs-C.
\end{equation}
A pair satisfying $\mathcal I_{\theta\eta}=0<g_{\theta\eta}$
forces $\Omega(T)$ regret. If $D_\A(s)\asymp s^\rho$, $\rho>1$,
the lower rate is $\Omega(T^{1-1/\rho})$.
\end{proposition}
\begin{proof}
Set $N=\lfloor T/2\rfloor$ and choose $\theta,\eta$ attaining
$D_\A(s)$, with $g_{\theta\eta}\ge s$.
Either $\Reg_T^\pi(\theta)\ge Ts$, or Lemma~\ref{lem:KL}
bounds the KL divergences up to length $N$ by
$C\{T(d_{\theta\eta}+s b_{\theta\eta})+1\}\le C$.
In the latter case, sum Lemma~\ref{lem:common} over $2n\le t<N$
and use \eqref{eq:energy} to get
$\Reg_T^\pi(\theta)+\Reg_T^\pi(\eta)\ge cTg_{\theta\eta}-C$.
Both cases give $\max_{i\in\{\theta,\eta\}}\Reg_T^\pi(i)\ge cTs-C$.
Taking $\inf_\pi$ proves \eqref{eq:twopoint}.
Use a fixed pair with $\mathcal I_{\theta\eta}=0<g_{\theta\eta}$
for the linear bound, or take $s\asymp T^{-1/\rho}$ for the power bound.
\end{proof}
\begin{proposition}\label{prop:loglower}
If $g_{\theta\eta}>0$ for some $\theta,\eta\in\operatorname{int}\A$,
then $\Reg_T(\A)\ge c\log(eT)$ for sufficiently large $T$.
\end{proposition}
\begin{proof}
Analyticity gives a closed segment $\{\theta_s:s\in J\}\subset
\operatorname{int}\A$ on which an entry of $\mathcal H_{\theta_s,k}$,
$0\le k\le2n$, has derivative bounded away from zero.
Uniformly for $s,s+h\in J$,
\begin{equation*}
 g_{\theta_s\theta_{s+h}}\asymp h^2,\quad
 d_{\theta_s\theta_{s+h}}+b_{\theta_s\theta_{s+h}}\le Ch^2.
\end{equation*}
Choose an interval $I$ of positive length with $I+[0,\kappa_0]\subset J$
for small $\kappa_0>0$. For finite
$M_T=\sup_{\theta\in\A}\Reg_T^\pi(\theta)$, set
$N=\lfloor T/2\rfloor$ and $h_t=\kappa_0(t+M_T+1)^{-1/2}$.
Lemma~\ref{lem:KL} bounds the KL divergence between $\theta_s$
and $\theta_{s+h_t}$ uniformly for $s\in I$ and $2n\le t<N$.
Since $I+h_t\subset J$, integrate Lemma~\ref{lem:common} over $I$,
sum over $t$, and use \eqref{eq:energy} to obtain
$\sum_{2n\le t<N}(t+M_T+1)^{-1}
\le C\int_J\{\Reg_T^\pi(\theta_s)+1\}\,ds
\le C(M_T+1)$.
Unless $M_T\gtrsim\log(eT)$, the left side is logarithmic,
a contradiction.
\end{proof}

\begin{proof}[Proof of Theorem~\ref{thm:main}.]
By \eqref{eq:info} and the common covering set for every pair with
$b_{\theta\eta}=0$, case (ii) holds exactly when some set is
incompatible. In case (iii), all sets are compatible and
Theorem~\ref{thm:exponent}(ii) defines $\rho_*$.
Lemma~\ref{lem:block} and Proposition~\ref{prop:global} give the upper bounds.
Proposition~\ref{prop:powerlower} gives $\Omega(T)$ in case (ii).
In case (iii), apply Proposition~\ref{prop:loglower} if $\rho_*=1$,
and Proposition~\ref{prop:powerlower} to a set attaining $\rho_*$
if $\rho_*>1$. These bounds also hold on $\K$.
\end{proof}

\section{Examples}\label{sec:examples}
\subsection{Logarithmic regret without identifying the parameter sign}
Take independent unit Gaussian noises, identity cost weights, and
\[
h_{t+1}=w^h_{t+1},\quad
 z_{t+1}=\theta h_t+w^z_{t+1},\qquad r_{t+1}=z_t+u_t+w^r_{t+1},
\]
$y_t=\col(z_t+v_t^z,r_t+v_t^r)$, $\theta\in[-M,M]$, $M>0$.
Put $q=1+\theta^2$, $q'=1+\eta^2$, and $\ell=q/(q+1)$.
In stationarity $\E_\theta[z_t\mid\F_t]=\ell y_t^z$, so the prediction residual is
$(y_t^z,y_t^r-u_{t-1}-\ell y_{t-1}^z)^\top$.
Only its lag-zero diagonal entries and lag-one lower-left entry
change, by $q'-q$, $(1-\ell)^2(q'-q)$ and $(1-\ell)(q'-q)$,
respectively. Therefore
\begin{equation}\label{eq:logexample}
g_{\theta\eta}=\frac{(q'-q)^2}{4(q+1)^2(q'+1)^2},\qquad b_{\theta\eta}=0,\qquad d_{\theta\eta}=(q'-q)^2\{1+(1-\ell)^2+(1-\ell)^4\}.
\end{equation}
Since $q,q'\ge1$, $d_{\theta\eta}\ge64g_{\theta\eta}$, giving $\Theta(\log T)$ regret.
The stationary optimal law and the observation law it induces depend only on $\theta^2$.

\subsection{A cube-root example}\label{sec:cuberoot}
Let
\begin{equation}\label{eq:cuberoot_model}
z_{t+1}=w^z_{t+1},\quad r_{t+1}=az_t+bu_t+w^r_{t+1},\qquad y_t=\col(z_t+v_t^z,r_t+v_t^r),
\end{equation}
with identity noise covariances and cost weights, zero initial state, and
\[
 \K=\{(a,b):a^2+(b-1)^2\le1/64\}.
\]
For $\theta=(a,b)$ and $\eta=(a',b')$, set
$\kappa(a,b)=ab/(1+b^2)$ and write
$\kappa=\kappa(a,b)$, $\kappa'=\kappa(a',b')$.
The differences in input responses and optimal laws are
$b_{\theta\eta}=(b-b')^2$ and
$g_{\theta\eta}=(\kappa-\kappa')^2/4$, respectively. A uniform polynomial certificate gives
\begin{equation}\label{eq:cubecert}
 g_{\theta\eta}^{\,3}\le C\mathcal I_{\theta\eta}^{\,2}
 \qquad\text{throughout }\K^2.
\end{equation}
Along $(a,b)=(t,1)$ and $(a',b')=(t+t^2/2,1+t)$, the quantities satisfy
$g_{\theta\eta}\asymp t^4$, $\mathcal I_{\theta\eta}\asymp t^6$. Thus $\rho_*=3/2$ and
\begin{equation}\label{eq:cuberoot}
 \Reg_T(\K)=\Theta(T^{1/3}).
\end{equation}
Appendix~\ref{app:examples} gives the certificate.

\subsection{A family without a common LTI stabilizer}
For $\theta\in[-1,1]$, take
\[
A_\theta=\begin{pmatrix}0&2\theta\\2\theta&0\end{pmatrix},
 \quad B_\theta=\binom{\theta}{0},\quad C=(1,0),\qquad W=Q=\id_2,\quad V=R=1.
\]
Every model is stabilizable and detectable. The endpoint plants have transfers $\pm\zeta/(\zeta^2-4)$. The closed-loop characteristic polynomials under a common proper LTI controller $N/D$, with monic $D$, would be $(\zeta^2-4)D\mp\zeta N$. Their values at $2$ have product $-4N(2)^2\le0$, whereas real monic Schur polynomials are positive at $2$. The controller $u_t=-4\eta y_{t-1}$ stabilizes the physical plant whenever $|\theta-\eta|\le1/8$.

\section{Conclusion}\label{sec:conclusion}
Minimax regret on $\K$ is bounded, logarithmic, linear, or
$\Theta(T^{1-1/\rho_*})$ with rational $\rho_*>1$.
Optimal control laws, residual covariances and input responses determine
this classification. Causal policies attain these rates without a supplied
common stabilizer. For real algebraic data, the regret class and exponent
can be determined exactly before observations are collected.

\appendix
\section{Finite reconstruction, transport and the geometric exponent}\label{app:geometry}
Constants are uniform over model indices in each $\A$ satisfying
conditions~\ref{cond:cover-geometry}--\ref{cond:cover-stability}.

\subsection{Residual responses and finite reconstruction}
Allow an arbitrary residual input $r_\nu=u+K_\nu\bar x^\nu$ in
\eqref{eq:implement}. The same physical realization gives
\begin{equation}\label{eq:residual_response}
{\widetilde y}_\nu=\mathcal T_{i,\nu}r_\nu+n_i^\nu,\qquad \mathcal T_{i,\nu}=H_{i,\nu}(\zeta \id-M_{i,\nu})^{-1}E_{e,i,\nu}.
\end{equation}
Here $n_i^\nu$ is the Gaussian response when $r_\nu=0$. Write
$S_i^\nu$ for its stationary spectrum. For $i=\nu$, prediction
cancels the input, so $\mathcal T_{\nu,\nu}=0$.

Each plant and controller has $n$ states. The controller state matrix
is $F_{c,i}(\id-J_iC_i)$ and need not be Schur. Each transfer difference in \eqref{eq:g-b} has order at most $2n$.
Cayley--Hamilton reconstructs its later coefficients from the given prefixes. Compactness
bounds every fixed extension of those prefixes.
Each $\mathcal T_{i,\nu}$ has the uniformly Schur $2n$-state
matrix \eqref{eq:phys_matrix}. The triangular map in
\eqref{eq:short_probe_identity} therefore bounds its difference
through lag $4n$ by the differences between the first $2n$ plant impulse response coefficients. A recurrence of order at most $4n$
for these transfer differences bounds the remaining Wiener tail. The
reverse bound follows from \eqref{eq:short_b}. Likewise, covariance
differences have order at most $4n$ at positive lags. Thus
\begin{equation}\label{eq:finite_transfer}
\|\mathcal T_{\theta,\nu}-\mathcal T_{\eta,\nu}\|_{\Wien}^2
  \asymp b_{\theta\eta},\qquad \|S_\theta^\nu-S_\eta^\nu\|_{\Wien}^2
  \asymp d^\nu_{\theta\eta}.
\end{equation}
The recurrence bounds use the common Schur margin of the compact
companion matrices.

\begin{lemma}\label{lem:transport}
For some $C_{\rm tr}>0$ and all $\theta,\eta,\nu\in\A$,
\begin{equation}\label{eq:transport2}
 d_{\theta\eta}\le C_{\rm tr}\{d^\nu_{\theta\eta}
                  +b_{\theta\eta}g_{\theta\nu}\}.
\end{equation}
\end{lemma}
\begin{proof}
Under $\pi^\theta$, \eqref{eq:twoport} gives
${\widetilde y}_\nu=\mathsf A_{\nu\theta}{\widetilde y}_\theta$ and
$r_\nu=\mathsf C_{\nu\theta}{\widetilde y}_\theta$. For $i=\theta,\eta$,
\eqref{eq:residual_response} therefore yields
\[
 {\widetilde y}_\theta=\mathcal L_i n_i^\nu,\qquad
 \mathcal L_i=(\mathsf A_{\nu\theta}
       -\mathcal T_{i,\nu}\mathsf C_{\nu\theta})^{-1}.
\]
Also ${\widetilde y}_i=\mathsf A_{i\theta}{\widetilde y}_\theta$ and
$n_i^\nu=\mathsf A_{i\nu}^{-1}{\widetilde y}_i$, so
$\mathcal L_i=\mathsf A_{i\theta}^{-1}\mathsf A_{i\nu}$ is uniformly stable.
Subtracting the inverse identities gives
\[
 \mathcal L_\theta-\mathcal L_\eta=\mathcal L_\theta
 (\mathcal T_{\theta,\nu}-\mathcal T_{\eta,\nu})
 \mathsf C_{\nu\theta}\mathcal L_\eta.
\]
By \eqref{eq:finite_transfer} and \eqref{eq:twoport-bounds},
$\|\mathcal L_\theta-\mathcal L_\eta\|_{\Wien}\le C\sqrt{b_{\theta\eta}g_{\theta\nu}}$.
Since $S_i^\theta=\mathcal L_iS_i^\nu \mathcal L_i^*$, subtracting these spectral
products gives
$\|S_\theta^\theta-S_\eta^\theta\|_{\Wien}
\le C(\sqrt{d^\nu_{\theta\eta}}+\sqrt{b_{\theta\eta}g_{\theta\nu}})$.
Square and apply \eqref{eq:finite_transfer} with $\nu=\theta$.
\end{proof}

\subsection{Rational error exponent}
Compatibility makes $\mathcal I_{\theta\eta}>0$ whenever $g_{\theta\eta}>0$. Polynomial identities and stabilizing Lyapunov inequalities give semialgebraic maps for finite response coefficients, Riccati solutions and covariance matrices. Compact minimization makes $D_\A(s)$ in \eqref{eq:profile} positive and semialgebraic for $s>0$.

For a fixed pair, write $g=g_{\theta\eta}$, $d=d_{\theta\eta}$, $b=b_{\theta\eta}$ and $\mathcal I=\mathcal I_{\theta\eta}$. For $g\ge s$, nonnegativity gives $d+sb\ge(s/g)(d+gb)$.
Thus $\mathcal I\ge cg^r$ implies $D_\A(s)\ge cs^r$.
Conversely, at $s=g$, $d+gb\ge D_\A(g)\ge cg^r$. Hence
\begin{equation}\label{eq:profileiff}
 D_\A(s)\gtrsim s^r\quad\Longleftrightarrow\quad
 \mathcal I_{\theta\eta}\gtrsim g_{\theta\eta}^{\,r}.
\end{equation}
Along a suitable short line segment in $\operatorname{int}\A$, nonconstancy gives
$g\asymp h^2$ and $d+b=O(h^2)$, hence $D_\A(s)\lesssim s$.
The positive one-variable semialgebraic function $D_\A$
has a rational leading power at the origin \citep[Lemma~2.3]{pham}. Thus
$D_\A(s)\asymp s^\rho$ for rational $\rho\ge1$,
identified by \eqref{eq:profileiff} with \eqref{eq:rho}.

\begin{proof}[Proof of Theorem~\ref{thm:exponent}.]
If the optimal law were constant on each $\A_i$, agreement on
intersections and connectedness of $\K$ would make it constant on $\K$,
since the cover is finite and its members are compact.
Thus, in the stated compatible case with $g_{\theta\eta}>0$ for some
$\theta,\eta\in\K$, the maximum defining $\rho_*$ is over a nonempty set.
When $\rho_*$ is defined, let $\{\mathcal B_j\}_{j=1}^N\in\mathfrak C(\K)$ be another
compatible cover with threshold $\widetilde\delta>0$,
and let $r$ be its maximum exponent.
Take $G\ge1$ bounding $g$ on $\K\times\K$.
Since $g\le G$, \eqref{eq:rho} and finiteness give a common $c>0$
with $\mathcal I_{\theta\eta}\ge c g_{\theta\eta}^r$ on every
$\mathcal B_j$, including sets with $g\equiv0$.
Pairs in $\A_i$ with $b_{\theta\eta}\le\widetilde\delta$ share
some $\mathcal B_j$. For the other pairs, \eqref{eq:info} gives
\[
 \mathcal I_{\theta\eta}\ge g_{\theta\eta}b_{\theta\eta}
 \ge\widetilde\delta g_{\theta\eta}
 \ge\widetilde\delta G^{1-r}g_{\theta\eta}^r.
\]
Thus $\rho_*\le r$. Reversing the covers gives equality.
\end{proof}

\subsection{Endpoint transport refinement}\label{app:endpoint-refinement}
If $\rho=1$, fix $\lambda>0$ with $\mathcal I_{\theta\eta}\ge\lambda g_{\theta\eta}$
on $\A$ and a transport constant $C_{\rm tr}$ in \eqref{eq:transport2}.
Choose $0<\beta<\lambda/2$ with $\beta\le\lambda/(16C_{\rm tr})$,
then refine in $\psi$ so that $b_{\theta\eta}\le\beta$ on each refined set. Thus
\[
 g_{\theta\eta}\le\frac2\lambda d_{\theta\eta}
 \le\frac{2C_{\rm tr}}\lambda
 (d^\nu_{\theta\eta}+\beta g_{\theta\nu}).
\]
Since $2C_{\rm tr}\beta/\lambda\le1/8$, this proves \eqref{eq:endpoint}.
The refinement is fixed independently of $T$ and the unknown parameter.

\section{Block moments and cost}\label{app:moment}
With $E_x=(\id_n,0)$, the reset implies
$Z_{b_j}=E_x^\top x_{b_j}$. Projecting the iterated
\eqref{eq:physical_real} onto the physical state gives
\begin{equation}\label{eq:blockrec}
 x_{b_{j+1}}=E_xM_{\theta,\nu_j}^{N_j}E_x^\top x_{b_j}+\chi_j.
\end{equation}
By \eqref{eq:reset_contraction} the coefficient is at most
$1/2$. The centered stable convolution $\chi_j$ has bounded moments of each fixed order
from Gaussian noise and bounded probes. Minkowski gives
\[
 \|x_{b_{j+1}}\|_{L^r}\le\tfrac12\|x_{b_j}\|_{L^r}+C_r.
\]
Stability with the controller fixed bounds moments within each block. Each reset has
initialization cost at most $C(1+\E\|x_b\|^2)$.

Choose $\nu_0$ fixed in each selected set and update from completed blocks.
Condition on past observations, inputs, probes and $x_b$.
Each scalar statistic in \eqref{eq:sample} is $\xi^\top A\xi$ plus
initialization terms, with fresh independent standard Gaussian and Rademacher
coordinates $\xi$ and $A=A^\top$. Schur stability gives
$\|A\|_F\le CN^{-1/2}$ and $\|A\|\le C/N$.
Integrating the tail bound in Theorem~1.1 of \citet{rudelsonvershynin2013} gives,
with conditional expectations and norms, for fixed $r\ge1$,
\begin{equation}\label{eq:quadratic_moment}
\|\xi^\top A\xi-\E(\xi^\top A\xi)\|_{L^{2r}}
 \le C_r(\|A\|_F+\|A\|) \le C_rN^{-1/2}.
\end{equation}
Average the $2r$th powers in \eqref{eq:quadratic_moment}.
Stability and the moment assumption in Lemma~\ref{lem:block} make initialization
terms $O(N^{-1})$ in $L^{2r}$. At lag $k$, zero padding omits $k\le4n$ products, giving factor
$1-k/N$ and bias $O(N^{-1})$. Fixed feature dimension and Jensen give
\eqref{eq:sample_moment} for all fixed $q>0$.

On a block with a fixed controller, $r_\nu=\sigma\,e$ and \eqref{eq:twoport} give
\[
 r_\theta=\mathsf C_{\theta\nu}{\widetilde y}_\nu+
          \sigma\,\mathsf D_{\theta\nu}e.
\]
Compactness and condition~\ref{cond:cover-stability} uniformly bound the spectrum of ${\widetilde y}_\nu$.
Thus \eqref{eq:twoport-bounds} and $\|a+b\|^2\le2\|a\|^2+2\|b\|^2$ give
stationary residual energy $O(g_{\theta\nu}+\sigma^2)$. Finite-horizon square completion
\eqref{eq:bellman}, uniform Riccati convergence and the uniformly
bounded prescribed initial-state moments yield \eqref{eq:blockcost} and
$V_T^*(\theta)=TJ_\theta^*+O(1)$ uniformly.

\section{Two-model realization and causal Gaussian likelihood}\label{app:KL}
\subsection{Control residual energy}\label{app:energy}
Let $\F_{-1}=\sigma(\omega)$,
$m_t^\theta=\E_\theta^\pi[x_t\mid\F_{t-1}]$ and
$m_0^\theta=\mu_\theta$. With Kalman gain $J_{\theta,t}$ initialized
at prediction covariance $\Sigma_\theta$, set
$L_{\theta,t}=A_\theta J_{\theta,t}$,
$\epsilon_t^\theta=y_t-C_\theta m_t^\theta$ and
$\widehat x_t^\theta=m_t^\theta+J_{\theta,t}\epsilon_t^\theta$.
For $\Delta_t=q_t^\theta-m_t^\theta$, the common inputs cancel:
\[
\Delta_{t+1}=(A_\theta-L_\theta C_\theta)\Delta_t
       +(L_\theta-L_{\theta,t})\epsilon_t^\theta,\qquad \Delta_0=-\mu_\theta.
\]
The steady estimate $\bar x_t^\theta=q_t^\theta+J_\theta{\widetilde y}_{\theta,t}$ satisfies
\[
 \bar x_t^\theta-\widehat x_t^\theta
 =(\id-J_\theta C_\theta)\Delta_t
       +(J_\theta-J_{\theta,t})\epsilon_t^\theta.
\]
Geometric convergence of the filter gains, bounded innovation
covariances and the Schur matrix $A_\theta-L_\theta C_\theta$ give uniformly
\[
 \sum_{t\ge0}\E_\theta^\pi
 \|\bar x_t^\theta-\widehat x_t^\theta\|^2\le C.
\]

For $a_t=u_t+K_{\theta,T-t}\widehat x_t^\theta$,
\eqref{eq:bellman} and $S_{\theta,h}\succeq R\succ0$ give
$\sum_{t<T}\E\|a_t\|^2\le C\Reg_T^\pi(\theta)$. The exact predictor obeys
\[
m_{t+1}^\theta
 =(A_\theta-B_\theta K_{\theta,T-t})m_t^\theta+B_\theta a_t +(A_\theta-B_\theta K_{\theta,T-t})
             J_{\theta,t}\epsilon_t^\theta.
\]
For $t<T/2$ and large $T$, its transitions are uniformly close
to the Schur matrix $A_\theta-B_\theta K_\theta$, so
\[
 \sum_{t<T/2}\E\|\widehat x_t^\theta\|^2
 \le C\{T+\Reg_T^\pi(\theta)+1\}.
\]
Control Riccati convergence gives
$\sup_{t<T/2}\|K_\theta-K_{\theta,T-t}\|\le Cr^{T/2}$ for some $0<r<1$.
The decomposition
\[
 r_{\theta,t}=a_t+(K_\theta-K_{\theta,T-t})\widehat x_t^\theta
                  +K_\theta(\bar x_t^\theta-\widehat x_t^\theta)
\]
then bounds the residual energy by
$C\Reg_T^\pi(\theta)+Cr^T\{T+\Reg_T^\pi(\theta)+1\}+C$.
Increasing $C$ covers bounded $T$ and proves \eqref{eq:energy}.

\subsection{Likelihood of the observed history}
Substitution in \eqref{eq:implement} with
$s_t=\col(q_t^\theta,q_t^\eta)$ yields \eqref{eq:twoport}:
\begin{equation}\label{eq:twoport-realization}
s_{t+1}=\mathsf M s_t+H_{\widetilde y}{\widetilde y}_{\theta,t}+H_rr_{\theta,t},\qquad {\widetilde y}_{\eta,t}={\widetilde y}_{\theta,t}+E_{\widetilde y} s_t,\qquad r_{\eta,t}=r_{\theta,t}+D_{\widetilde y}{\widetilde y}_{\theta,t}+E_r s_t,
\end{equation}
with
\[
\begin{gathered}
\mathsf M=\begin{pmatrix}F_{c,\theta}&0\\L_\eta C_\theta-B_\eta K_\theta&F_{o,\eta}\end{pmatrix},\qquad H_{\widetilde y}=\binom{F_{c,\theta}J_\theta}{L_\eta-B_\eta K_\theta J_\theta},
\quad H_r=\binom{B_\theta}{B_\eta},\\
E_{\widetilde y}=(C_\theta,-C_\eta),\quad D_{\widetilde y}=K_\eta J_\eta-K_\theta J_\theta,\qquad E_r=(-K_\theta+K_\eta J_\eta C_\theta,
 K_\eta(\id-J_\eta C_\eta)).
\end{gathered}
\]
Both diagonal blocks of $\mathsf M$ are uniformly Schur, and
exchanging models inverts the map. Eliminating $q^i$ yields
${\widetilde y}_i=\mathcal R_i(y-\mathcal G_i u)$ and
$r_i=\mathcal U_i(u-\mathcal H_i y)$, where
\[
\mathcal R_i=\id-C_i(\zeta\id-F_{o,i})^{-1}L_i,\qquad \mathcal U_i=\id+K_i(\id-J_iC_i)(\zeta\id-F_{o,i})^{-1}B_i.
\]
The inverse in \eqref{eq:implement} gives the maps
$y=\mathcal Y_i{\widetilde y}_i$ when $r_i=0$ and
$u=\mathcal D_i r_i$ when ${\widetilde y}_i=0$, with
\[
\mathcal Y_i=\id+C_i(\zeta\id-F_{c,i})^{-1}F_{c,i}J_i,\qquad \mathcal D_i=\id-K_i(\zeta\id-F_{c,i})^{-1}B_i.
\]
Consequently,
\[
\mathsf C_{\eta\theta}=\mathcal U_\eta(\mathcal H_\theta-\mathcal H_\eta)\mathcal Y_\theta,\qquad \mathsf B_{\eta\theta}=\mathcal R_\eta(\mathcal G_\theta-\mathcal G_\eta)\mathcal D_\theta.
\]
The maps with identity diagonal blocks compare the squared prefix norms of
$\mathsf B_{\eta\theta,1:2n},\mathsf C_{\eta\theta,0:2n}$ with $b_{\theta\eta},g_{\theta\eta}$. Their uniformly
Schur $2n$-state recurrences bound the remaining Wiener tails by
these prefixes. Stability of $\mathsf D_{\eta\theta}$ gives its coefficient
decay. This proves \eqref{eq:twoport-bounds}.

For $i,\nu\in\A$, set $r_\nu=0$ and drive the
matched predictor by ${\widetilde y}_i$. In state $(q^i,q^\nu)$ the transition
is $M_{i,\nu}$ from \eqref{eq:phys_matrix}, and
${\widetilde y}_\nu=H_{i,\nu}\col(q^i,q^\nu)+{\widetilde y}_i$. Consequently
\begin{equation}\label{eq:innovation_inverse}
 \mathsf A_{i\nu}^{-1}=\id+H_{i,\nu}(\zeta \id-M_{i,\nu})^{-1}
 \binom{L_i-B_iK_\nu J_\nu}{F_{c,\nu}J_\nu}.
\end{equation}
The uniform Schur bound for $M_{i,\nu}$ proves stability of $\mathsf A_{i\nu}^{-1}$.

To compare the residual covariances, write
\[
 {\widetilde y}_{i,0:N-1}=\col({\widetilde y}_{i,0},\ldots,{\widetilde y}_{i,N-1}).
\]
Under true model $i$, the predictor error
$\widetilde x_{i,t}=x_t-q_t^i$ satisfies
\[
 \widetilde x_{i,t+1}=F_{o,i}\widetilde x_{i,t}+w_{t+1}-L_iv_t,
 \qquad {\widetilde y}_{i,t}=C_i\widetilde x_{i,t}+v_t.
\]
These processes are policy-independent. With zero initial error,
the noise-to-residual matrix $F_N$ has right inverse
$z_t\mapsto(L_iz_t,z_t)$, which keeps $\widetilde x_{i,t}=0$ and outputs $z_t$.
Its norm is $\sqrt{1+\|L_i\|^2}$, so
$F_NF_N^\top\succeq(1+\|L_i\|^2)^{-1}\id$.
The noise lower bounds and independent initial error give the lower
bound below. Schur stability and bounded initial moments give the
upper bound and summable mean energy:
\begin{equation}\label{eq:direct-covbounds}
 c\id\preceq\Cov_i({\widetilde y}_{i,0:N-1})\preceq C\id.
\end{equation}
The same right inverse for stationary convolutions gives
$c\id\preceq S_i^i\preceq C\id$. Under $\pi^\nu$,
${\widetilde y}_i=\mathsf A_{i\nu}{\widetilde y}_\nu$. Stability of $\mathsf A_{i\nu}$ and of its inverse in
\eqref{eq:innovation_inverse} transfers the finite covariance bounds
and gives $S_i^\nu=\mathsf A_{i\nu}^{-1}S_i^i(\mathsf A_{i\nu}^{-1})^*$
with the same spectral bounds.
Under the frozen controller $\pi^\nu$, define
\[
 \Sigma_{i,N}^{(\nu)}=\Cov_i^{\pi^\nu}({\widetilde y}_{\nu,0:N-1})
 \in\mathbb R^{Np\times Np}.
\]
Summability of the observability columns of the $2n$-state realization
bounds the difference between
the covariance under the specified initial distribution and the
stationary covariance. The recurrence of order at most $4n$ for the covariance differences
and \eqref{eq:d} then give
\begin{equation*}
\norm{\Sigma_{\theta,N}^{(\theta)}-
\Sigma_{\eta,N}^{(\theta)}}_F^2\le C(Nd_{\theta\eta}+1).
\end{equation*}
Gaussian KL and the mean and covariance bounds yield
\begin{equation}\label{eq:frozen-KL}
D_{\mathrm{KL},N}(\theta,\eta)\le C(Nd_{\theta\eta}+1).
\end{equation}

For an arbitrary measurable policy, fix $N$ and write ${\widetilde y}_i={\widetilde y}_{i,0:N-1}$ and
$m_i=\E_i{\widetilde y}_i\in\mathbb R^{Np}$. These Gaussian means are
policy-independent. Let $T_i$ be the lower triangular Gaussian
prediction matrix with identity diagonal blocks such that
$\epsilon_i=T_i({\widetilde y}_i-m_i)$ has independent blocks, and set
$\Omega_{i,t}=\Cov_i(\epsilon_{i,t})$ for $0\le t<N$.
Schur complements of \eqref{eq:direct-covbounds} give
$c\id\preceq\Omega_{i,t}\preceq C\id$. The identity
$T_i\Cov({\widetilde y}_i)T_i^\top=\diag_t\Omega_{i,t}$ gives $\|T_i\|\le C$.
For the following formulas, use the lower triangular convolution matrices
of the transfer operators on $0,\ldots,N-1$ and write
$r_\theta=\col(r_{\theta,0},\ldots,r_{\theta,N-1})$.
On every observed history, \eqref{eq:twoport} gives
\[
\epsilon_\eta-\epsilon_\theta=a+b,\qquad a=(T_\eta\mathsf A_{\eta\theta}-T_\theta){\widetilde y}_\theta
      +T_\theta m_\theta-T_\eta m_\eta,\qquad b=T_\eta\mathsf B_{\eta\theta} r_\theta.
\]
Both $a_t,b_t$ are predictable because the operator differences
have zero diagonal. We compare the two models at the same value of
$\omega$. The triangular maps with identity diagonal blocks then preserve Gaussian
conditional likelihoods for arbitrary measurable actions.
With $\norm{z}_\eta^2=\sum_{t<N}z_t^\top\Omega_{\eta,t}^{-1}z_t$,
\[
\mathrm{KL}(\Prob_{\theta,N}^\pi\Vert\Prob_{\eta,N}^\pi)
=\mathcal D_N+\tfrac12\E_\theta\norm{a+b}_\eta^2,
\]
where the centered Gaussian divergence is
\[
 \mathcal D_N=\sum_{t<N}\mathrm{KL}
 \bigl(N(0,\Omega_{\theta,t})\Vert N(0,\Omega_{\eta,t})\bigr).
\]
Under $\theta$, $a$ is policy-independent, and
$D_{\mathrm{KL},N}(\theta,\eta)=\mathcal D_N+
\frac12\E_\theta\norm{a}_\eta^2$ for the frozen controller.
Hence $\|a+b\|_\eta^2\le2\|a\|_\eta^2+2\|b\|_\eta^2$, \eqref{eq:twoport-bounds}
and \eqref{eq:frozen-KL} yield
\[
\mathrm{KL}\le2D_{\mathrm{KL},N}(\theta,\eta)+
 Cb_{\theta\eta}\E_\theta\sum_{t<N}\norm{r_{\theta,t}}^2,
\]
which proves Lemma~\ref{lem:KL} using \eqref{eq:energy}.

\section{Cover construction and supervisory details}\label{app:cover}
\subsection{Stability for models with the same response}
Models with the same response vector $\psi$ have identical
transfers from inputs to outputs by finite reconstruction. Their minimal
realizations are similar. The remaining modes in the
Kalman decomposition are Schur by stabilizability and detectability.
These models therefore share stabilizers. Compactness and eigenvalue
continuity preserve stability in response neighborhoods with a uniform
positive radius.

Take a finite $\varepsilon$-net $\{\psi(a_i)\}$ of $\psi(\K)$,
with $a_i\in\K$ and $8\varepsilon$ below the stability radius.
Initially set $\varepsilon_i=\varepsilon$ in
\begin{equation}\label{eq:cells}
 \A_i=\overline{\{\theta\in\operatorname{int}\K:
 \norm{\psi(\theta)-\psi(a_i)}<4\varepsilon_i\}}.
\end{equation}
These sets are compact, semialgebraic and closures of their interiors.
For compatible initial sets with exponent $1$, choose
$0<\tilde\varepsilon<\varepsilon/2$ ensuring \eqref{eq:endpoint}
on subsets with response diameter $\le8\tilde\varepsilon$
(Appendix~\ref{app:geometry}). Replace each such set, centered at $a$,
by sets \eqref{eq:cells} of radius $4\tilde\varepsilon$ centered on a
finite $\tilde\varepsilon$-net of
$\{\psi(\theta):\theta\in\K,\ \norm{\psi(\theta)-\psi(a)}\le2\varepsilon\}$, with centers in this region.
New sets lie in their parent since
$2\varepsilon+4\tilde\varepsilon<4\varepsilon$ and satisfy
\eqref{eq:endpoint}, having response diameter $\le8\tilde\varepsilon$.
Retain the other sets and relabel as $a_i,\varepsilon_i,\A_i$.
Every $\theta\in\K$ has a final center with
$\norm{\psi(\theta)-\psi(a_i)}\le\varepsilon_i$.
If $b_{\theta\eta}\le\varepsilon_0^2$, then
$\norm{\psi(\eta)-\psi(a_i)}\le2\varepsilon_i<4\varepsilon_i$.
Continuity and $\K=\overline{\operatorname{int}\K}$ give
$\theta,\eta\in\A_i$, proving condition~\ref{cond:cover-close}
with $\delta=\varepsilon_0^2$.

\subsection{Pulse response estimation after each restart}
Fix the true model and construct the global stable observer $\check x_0=0$,
\[
 \check x_{t+1}=(A_\theta-L_\theta C_\theta)\check x_t
 +B_\theta u_t+L_\theta y_t.
\]
Thus $\norm{\check x_t}\le C R_t$ for the same running bound
used in the monitor. Set $\check y_t=y_t-C_\theta\check x_t$. The stable observer-error recursion and Cayley--Hamilton bound its
unknown initial observable contribution by the first $n$ residuals
and the noise envelope. Thus for fixed $0<r<1$,
\begin{equation}\label{eq:observerbound}
\norm{\check y_t}\le C(R_n+D_N),\qquad \norm{x_t}\le C\{R_t+D_N+r^t\norm{x_0}\},\qquad t\le N.
\end{equation}
Also,
\[
 \check x_{t+1}=A_\theta\check x_t+B_\theta u_t+L_\theta\check y_t,
 \quad y_t=C_\theta\check x_t+\check y_t.
\]
At each restart $\|\check x_{\tau_j}\|\le C\overline R_j$.
Over the fixed $m\ell+1$ pulse steps, $A_\theta^k$ is uniformly
bounded. Dividing the $k$-step response to pulse $r$ in
\eqref{eq:pulses} by $C_{\rm p}\Lambda_ja_{\rm p}^r$ yields
$C_\theta A_\theta^{k-1}B_\theta$ with error at most
\begin{equation}\label{eq:trialpulseerror}
 C_{\rm est}\left\{\frac{\overline R_j+D_N}{C_{\rm p}(\overline R_j+h_j)}+\frac1{a_{\rm p}}\right\},
 \qquad 1\le k\le\ell.
\end{equation}
In \eqref{eq:trialpulseerror}, $C_{\rm est}>0$ depends only on the known family. The terms bound contributions
from the initial observer state and noise, and from earlier pulses, respectively.
First choose $C_{\rm p},a_{\rm p}$ so that the estimated response vector
has error $<\varepsilon_0$ when $D_N\le h_j$.
Some center satisfies
$\norm{\widehat\psi-\psi(a_i)}/\varepsilon_i<2$.
The selected index therefore satisfies
$\norm{\psi(\theta)-\psi(a_i)}
<2\varepsilon_i+\varepsilon_0\le3\varepsilon_i<4\varepsilon_i$,
so $\theta\in\A_i$. For these fixed pulses, pulse and block stability give
a known constant $C_{\rm cont}>0$ bounding input and output norms by
$C_{\rm cont}(\Lambda_j+D_N)\le2C_{\rm cont}\Lambda_j$ when $D_N\le h_j$.
Fix $C_{\rm out}>2C_{\rm cont}$ to prevent another trigger.

\subsection{Restart count and rare-event cost}
Before rejection the monitor bounds inputs by
$C_{\rm out}\Lambda_j$, even on a wrong set.
Equation~\eqref{eq:observerbound} bounds the next output overshoot.
The fixed pulse bound controls the next trial. Hence
\begin{equation}\label{eq:Rrec}
 \overline R_{j+1}\le A(\overline R_j+h_j+D_N),\qquad A\ge2.
\end{equation}
Let $L_N=\lceil\log_2^+(D_N/h_0)\rceil$. If $j\ge L_N$,
then $h_j\ge D_N$ and no rejection occurs. Hence the last
trial begun before $N$ has index at most $L_N$. With
$E_j=\overline R_j+h_j+D_N$, \eqref{eq:Rrec} gives
\[
 E_{j+1}\le(A+2)E_j,\qquad E_0\le2h_0+D_N.
\]
Thus $q\ge1+\log_2(A+2)$ gives \eqref{eq:pathbound}.
\eqref{eq:observerbound} bounds state energy by a fixed
polynomial in $h_0+D_N$ plus $r^{2t}\|x_0\|^2$.
The tail bound \eqref{eq:noise_tail} and Cauchy--Schwarz make
the exceptional-event cost in Proposition~\ref{prop:global} $O(1)$.

\section{Optimal control laws under output costs}\label{app:output}
Expected cumulative output cost equals the expected cumulative state cost with
weight $C_i^\top Q_yC_i$ plus $T\tr(Q_yV_i)$.
Detectability selects the stabilizing Riccati branch even if this
state weight is singular. Fix a covering set $\A$ and $\theta,\eta\in\A$.
Set $(z,\delta)=({\widetilde y}_\theta,r_\theta)$ and write the stable inverse
of \eqref{eq:implement} as
\[
 \binom{Q_y^{1/2}y}{R^{1/2}u}=Fz+G\delta.
\]
In particular, $\mathcal T_{\theta,\theta}=0$.
For $i=\theta,\eta$, setting $\delta=Hn_i^\theta$, $z=(\id+\mathcal T_{i,\theta}H)n_i^\theta$ gives, with $\mathcal B_i=F\mathcal T_{i,\theta}+G$,
\begin{equation*}
 \mathcal J_i(H)=\|F+\mathcal B_iH\|_{S_i^\theta}^2,\qquad \mathcal B_\theta=G.
\end{equation*}
Here $\|U\|_S^2=\langle U,U\rangle_S$, where for transfer maps $U,V$,
\[
 \langle U,V\rangle_S=(2\pi)^{-1}\operatorname{Re}
 \int_{-\pi}^\pi\tr(USV^*)\,d\omega.
\]
All maps in this integral are evaluated at $\zeta=e^{\mathrm i\omega}$.
The stable map $\mathcal L$ recovering $\delta$ from the weighted outputs
has $\mathcal LF=0$, $\mathcal LG=\id$, hence $\mathcal L\mathcal B_i=\id$.
Together with $S_i^\theta\succeq c\id$, this gives
\begin{equation}\label{eq:coercive}
 \|\mathcal B_iH\|_{S_i^\theta}^2\ge c\|H\|_{H_2}^2.
\end{equation}
At $\theta$, $n_\theta^\theta=z$, so $\delta=\epsilon Hz$ is feasible
for every causal $H$ with a finite impulse response. Optimality of zero gives
\[
 0\le \mathcal J_\theta(\epsilon H)-\mathcal J_\theta(0)
 =2\epsilon\langle F,GH\rangle_{S_\theta^\theta}
  +\epsilon^2\|GH\|_{S_\theta^\theta}^2.
\]
Both signs of $\epsilon$ force $\langle F,GH\rangle_{S_\theta^\theta}=0$.
Truncation extends this identity to every causal $H\in H_2$.
Under $\pi^\eta$, $\delta=\mathsf C_{\theta\eta}{\widetilde y}_\eta$ and
$n_\eta^\theta=\mathsf A_{\eta\theta}^{-1}{\widetilde y}_\eta$, so its stable response is
$H=\mathsf C_{\theta\eta}\mathsf A_{\eta\theta}$.
Expanding $\mathcal J_\eta(H)\le \mathcal J_\eta(0)$ and using \eqref{eq:coercive} gives
\[
\begin{gathered}
c\|H\|_{H_2}^2
 \le-2\langle F,\mathcal B_\eta H\rangle_{S_\eta^\theta} =-2\langle F,F\mathcal T_{\eta,\theta} H\rangle_{S_\eta^\theta}
   -2\langle F,GH\rangle_{S_\eta^\theta-S_\theta^\theta}\\
\le C(\|\mathcal T_{\eta,\theta}\|_{\Wien}+\|S_\eta^\theta-S_\theta^\theta\|_{\Wien})\|H\|_{H_2}.
\end{gathered}
\]
The last step uses Cauchy--Schwarz and uniform stable bounds.
Since $\mathsf C_{\theta\eta}=H\mathsf A_{\eta\theta}^{-1}$,
\eqref{eq:twoport-bounds} gives $g_{\theta\eta}\le C\|H\|_{H_2}^2$.
Equation~\eqref{eq:finite_transfer} gives
$\|S_\eta^\theta-S_\theta^\theta\|_{\Wien}^2\asymp d_{\theta\eta}$ and
$\|\mathcal T_{\eta,\theta}\|_{\Wien}^2\asymp b_{\theta\eta}$.
Divide when $H\ne0$ (the zero case is immediate) and square to obtain
\begin{equation}\label{eq:outputstab}
 g_{\theta\eta}\le C(d_{\theta\eta}+b_{\theta\eta}).
\end{equation}
Thus $d_{\theta\eta}=b_{\theta\eta}=0$ implies $g_{\theta\eta}=0$.
Moreover,
\[
\mathcal I_{\theta\eta}
 \ge\min(1,g_{\theta\eta})(d_{\theta\eta}+b_{\theta\eta}) \ge g_{\theta\eta}^2/[C\max(1,G_\A)].
\]
This proves Corollary~\ref{cor:output}.

\section{Cube-root certificate and exact algebraic classification}\label{app:examples}
For \eqref{eq:cuberoot_model}, abbreviate $g=g_{\theta\eta}$,
$d=d_{\theta\eta}$ and $\mathcal I=\mathcal I_{\theta\eta}$. Set
\[
 h=b'-b,\qquad e=a'-a-\kappa h,\qquad
 \tau=(a'^2-a^2+e^2)/2.
\]
In stationarity $\E_\theta[z_t\mid\F_t]=y_t^z/2$, so the prediction residual is
$(y_t^z,y_t^r-a y_{t-1}^z/2-bu_{t-1})^\top$.
Under $\eta$ and $u_t=-\kappa y_t^z/2$, its second component is
\[
 (a'+e)z_{t-1}/2-(a'-e)v_{t-1}^z/2+w_t^r+v_t^r.
\]
Its variance is $2+(a'^2+e^2)/2$ and its lag-one covariance
with $y^z$ is $e$. All other differing entries vanish. Hence
\begin{equation}\label{eq:explicit_cube}
g=\tfrac14(\kappa'-\kappa)^2,\quad b_{\theta\eta}=h^2,\qquad d=e^2+\tau^2,\quad \mathcal I=d+gh^2.
\end{equation}
The identity
\begin{equation}\label{eq:kappa_identity}
 \kappa'-\kappa
 =\frac{b'}{1+b'^2}e
 +\frac{a h}{(1+b^2)(1+b'^2)}
\end{equation}
is used with $f=b/(1+b^2)$ and $z=ah$.
On $|a|,|a'|,|b-1|,|b'-1|\le1/8$, we have
$f\in[2/5,1/2]$ and $|h|\le1/4$. Since $a'=a+fz+e$,
\[
 2\tau h=f(2+fh)z^2+2(1+fh)ez+2he^2
 \ge3z^2/8-31e^2/8,
\]
using $9|ez|/4\le3z^2/8+27e^2/8$.
Thus $z^2\le(31e^2+16|\tau h|)/3$.
Since $(1+b^2)(1+b'^2)\ge(113/64)^2>3$, the coefficients in
\eqref{eq:kappa_identity} are at most $1/2$ and $1/3$. Hence
\[
 g\le e^2/8+z^2/18
 \le151e^2/216+8|\tau h|/27\le d+\sqrt{dh^2}.
\]
Since $|\kappa|,|\kappa'|\le1/16$, $g\le1$. Multiplying by
$\sqrt g$ and using $2\sqrt{dgh^2}\le d+gh^2$ gives
\begin{equation}\label{eq:cube_exact}
 g^{3/2}\le d+\sqrt{dgh^2}\le\tfrac32\mathcal I,
 \qquad g^3\le\tfrac94\mathcal I^2.
\end{equation}
For $0<t\le1/16$, take $(a,b)=(t,1)$, $(a',b')=(t+t^2/2,1+t)$:
\[
 e=0,\quad \tau=t^3(1+t/4)/2,\quad
 g=\frac{t^4}{16(2+2t+t^2)^2}.
\]
Thus $g\sim t^4/64$ and $\mathcal I\sim17t^6/64$.
Every $r<3/2$ fails \eqref{eq:rho}. \eqref{eq:cube_exact} proves equality.

\subsection{Exact algebraic classification}\label{app:compute}
Starting at $\varepsilon=1$, halve $\varepsilon$ until, for all
$\theta,\nu\in\K$ with $\|\psi(\theta)-\psi(\nu)\|\le8\varepsilon$,
the matrix $M_{\theta,\nu}$ in \eqref{eq:phys_matrix} satisfies
\[
 \exists H=H^\top\succ0:\quad
 H-M_{\theta,\nu}^\top H M_{\theta,\nu}\succ0.
\]
Stability for models with the same $\psi$ and compactness in Appendix~\ref{app:cover} give a
positive admissible radius, so the search terminates.
On this pair domain, search integers $\ell\ge1$ until every pair
admits a symmetric $H$ with
\[
\id\preceq H\preceq2^\ell\id,\qquad M_{\theta,\nu}^\top H M_{\theta,\nu}
   \preceq(1-2^{-\ell})^2H.
\]
Uniform Schur stability guarantees termination. Set
$C_0=2^{\ell/2}$ and $\gamma=1-2^{-\ell}$, giving
$\|M_{\theta,\nu}^k\|\le C_0\gamma^k$, and choose
an integer $N_0>8n+2$ with $C_0\gamma^{N_0}\le1/2$
as in Lemma~\ref{lem:block}.

Starting with an algebraic $a_1\in\K$, extend $a_1,\ldots,a_j$
by an algebraic witness $a_{j+1}=\theta\in\K$ to
$\|\psi(\theta)-\psi(a_i)\|\ge\varepsilon$ for every $i\le j$
while one exists.
Compact packing forces termination, giving a strict
$\varepsilon$-net and the initial sets in \eqref{eq:cells}
with $\varepsilon_i=\varepsilon$.

Test constancy of the optimal law by $g_{\theta\eta}=0$ for all
$\theta,\eta\in\K$, and test each set for incompatibility by
$\exists\theta,\eta\in\A:\ \mathcal I_{\theta\eta}=0,
g_{\theta\eta}>0$. In compatible sets with varying optimal laws, quantifier
elimination gives the graph $y=D_\A(s)$ in terms of finitely many polynomial sign conditions. The product $H(s,y)$ of the nonzero polynomials in these sign conditions vanishes along the graph. If $D_\A(s)\sim as^\rho$ with $a>0$, at least two monomials $c_{ij}s^iy^j$ of $H$ have equal smallest weight $i+j\rho$. Hence $\rho$ belongs to the finite candidate set
\[
 \left\{\frac{i'-i}{j-j'}:\ c_{ij}c_{i'j'}\ne0,
 j\ne j',\ \frac{i'-i}{j-j'}\ge1\right\}.
\]
Test each reduced candidate $p/q$ in increasing order:
\begin{equation}\label{eq:exp_certificate}
 \exists C>0\ \forall\theta,\eta\in\A:\qquad
 g_{\theta\eta}^{p}\le C\mathcal I_{\theta\eta}^{q}.
\end{equation}
The first passing value is $\rho$. When $\rho=1$, use universal
tests to select $\lambda>0$ with
$\mathcal I_{\theta\eta}\ge\lambda g_{\theta\eta}$ for all
$\theta,\eta\in\A$, and a transport constant $C_{\rm tr}$
in \eqref{eq:transport2}.
Choose $0<\beta<\lambda/2$ with $\beta\le\lambda/(16C_{\rm tr})$
on each such set. Take algebraic $0<\tilde\varepsilon<\varepsilon/2$
with $64\tilde\varepsilon^2\le\beta$ for every chosen $\beta$.
Apply the same witness search to the inner regions in
Appendix~\ref{app:cover}. New sets with varying laws have exponent $1$.

\paragraph{Acknowledgments.} This work was supported in part by the National Natural Science Foundation of China under Grant 12401585, the Guangdong Basic and Applied Basic Research Foundation under Grant 2024A1515011542, and the General Program of Shenzhen Natural Science Foundation under Grant JCYJ20250604181037012.

\begingroup
\interlinepenalty=10000

\endgroup
\end{document}